\documentclass{gtart_a}
\pdfoutput=1
\usepackage{pinlabel,amscd}

\title[Tightness and efficiency of automorphisms of
handlebodies]{Tightness and efficiency of irreducible\\automorphisms of
handlebodies}

\author{Leonardo Navarro Carvalho}
\givenname{Leonardo Navarro}
\surname{Carvalho}
\address{Departamento de Matem\'{a}tica Aplicada\\
IM---Universidade Federal Fluminense\\\newline
Niter\'{o}i, RJ\\Brazil}
\email{leonardo\_carvalho@vm.uff.br}
\urladdr{}

\volumenumber{10}
\issuenumber{}
\publicationyear{2006}
\papernumber{3}
\lognumber{0493}
\startpage{57}
\endpage{95}

\doi{}
\MR{}
\Zbl{}


\received{1 September 2004}
\revised{15 December 2005}
\accepted{26 November 2005}
\published{4 March 2006}
\publishedonline{4 March 2006}
\proposed{Cameron Gordon}
\seconded{Dave Gabai, Jean-Pierre Otal}
\corresponding{}
\editor{}
\version{}

\begin{asciiabstract}
Among (isotopy classes of) automorphisms of handlebodies those called
irreducible (or generic) are the most interesting, analogues of
pseudo-Anosov automorphisms of surfaces.  We consider the problem of
isotoping an irreducible automorphism so that it is most efficient
(has minimal growth rate) in its isotopy class. We describe a
property, called tightness, of certain invariant laminations, which we
conjecture characterizes this efficiency. We obtain partial results
towards proving the conjecture.  For example, we prove it for genus
two handlebodies.  We also show that tightness always implies
efficiency.

In addition, partly in order to provide counterexamples in our study
of properties of invariant laminations, we develop a method for
generating a class of irreducible automorphisms of handlebodies.
\end{asciiabstract}

\AtBeginDocument{\def\notin{\not\in}}

\makeatletter
\def\cnewtheorem#1[#2]#3{\newtheorem{#1}{#3}[section]
\expandafter\let\csname c@#1\endcsname\c@thm}

\makeatother  %  move after \newtheorem block

\newtheorem{theorem}{Theorem}
\newtheorem{lemma}[theorem]{Lemma}
\newtheorem{proposition}[theorem]{Proposition}
\newtheorem{corollary}[theorem]{Corollary}

\newtheorem{conjecture}[theorem]{Conjecture}
\newtheorem*{conjecture*}{Conjecture}

\newtheorem*{problem*}{Problem}
\newtheorem*{proposition*}{Proposition}
\newtheorem*{no label}{}
\newtheorem*{conj-equiv}{\fullref{Cj:equivalence}}
\newtheorem*{thm-minim}{\fullref{T:minimality}}
\newtheorem*{corr-true}{\fullref{C:true genus 2}}
\newtheorem*{corr-small}{\fullref{C:small growth}}
\newtheorem*{corr-growth}{\fullref{C:growth power}}
\newtheorem*{prop-valence}{\fullref{P:valence}}

\theoremstyle{definition}
\newtheorem{definition}[theorem]{Definition}
\newtheorem{example}[theorem]{Example}
\newtheorem{process}{Process}
\newtheorem{obs}[theorem]{Remark}
\newtheorem{obss}[theorem]{Remarks}

\numberwithin{theorem}{section}

\def\ring#1{\mathring#1}
\def\cal#1{\mathcal{#1}}

\def\tightness{tightness}
\def\tight{tight}
\def\tightening{tightening}

\def\inter#1#2#3{{#1}\bullet({#2},{#3})}

\def\overA{{}\mskip2mu\overline{\mskip-2mu A\mskip-2mu}\mskip2mu}

\begin{document}

\begin{abstract}
Among (isotopy classes of) automorphisms of handlebodies those
called irreducible (or generic) are the most interesting,
analogues of pseudo-Anosov automorphisms of surfaces.  We
consider the problem of isotoping an irreducible automorphism so
that it is most efficient (has minimal growth rate) in its isotopy
class. We describe a property, called tightness, of certain
invariant laminations, which we conjecture characterizes this
efficiency. We obtain partial results towards proving the
conjecture.   For example, we prove it for genus two handlebodies.
We also show that tightness always implies efficiency.

In addition, partly in order to provide counterexamples in our
study of properties of invariant laminations, we develop a method
for generating a class of irreducible automorphisms of
handlebodies.
\end{abstract}

\maketitle

\section{Introduction}\label{S:intro}

\subsection{Some history and background}\label{SS:background}
The classification of \emph{automorphisms} (ie
self-diffeomorphisms) of a manifold, up to isotopy, is a
fundamental problem. Nielsen addressed the case where the manifold
is a compact and connected surface and his results were later
substantially improved by Thurston (see the work of Nielsen
\cite{Nielsen:SurfacesI,Nielsen:SurfacesII,Nielsen:SurfacesIII},
Thurston \cite{WPT:SurfaceAutos}, Handel and Thurston \cite{MHWT:Nielsen}
and Fathi, Laudenbach and Poenaro \cite{FLP:TravauxThurston}).
We briefly state their main result: An automorphism of a surface
is, up to isotopy, either \emph{periodic} (ie has finite order),
\emph{reducible} (ie preserves an essential codimension 1
submanifold) or \emph{pseudo-Anosov}.  We refer the reader to any
of
\cite{FLP:TravauxThurston,MHWT:Nielsen,WPT:SurfaceAutos}, or the book by
Casson and Bleiler \cite{ACSB:SurfaceAutos}
for details, including the definition of a pseudo-Anosov
automorphism. The Nielsen--Thurston theory also shows that the
reducible case may be reduced to the other two cases. Since
periodic automorphisms are relatively easy to understand, the
remaining \emph{irreducible} case, ie the pseudo-Anosov case, is
the most interesting and rich one. Indeed, pseudo-Anosov
automorphisms of surfaces are the subject of intense and wide
research (see \cite{WPT:SurfaceAutos}).

The author and U Oertel have attempted a similar classification
of automorphisms of three-dimensional manifolds. First Oertel
\cite{UO:Autos}
studied the special case of automorphisms of handlebodies and
compression bodies. Recently, the author and
Oertel \cite{LCUO:classification} have described a classification of
automorphisms of an
arbitrary compact, connected, and orientable three-manifold $M$
satisfying Thurston's Geometrization Conjecture.
The classification falls short of
being a perfect analogue of the Nielsen--Thurston classification
because an automorphism of a reducible three-manifold must, in
general, first be written as a composition of two automorphisms,
each of which fits into the classification. With the help of
standard three-manifolds techniques (eg  the
JSJ-decomposition, Bonahon's ``characteristic compression body'',
Seifert fiberings, and Heegaard splittings), the authors describe
how a given automorphism of a three-manifold $M$ may be
decomposed along suitable ``reducing surfaces''. The automorphisms
which arise and were not previously understood are, in fact, just
automorphisms of \emph{handlebodies} and \emph{compression bodies}
(defined below). Therefore the path to a better understanding of
automorphisms of arbitrary compact three-manifolds leads directly
to the study of automorphisms of handlebodies and compression
bodies. The goal of this paper is to improve our understanding of
automorphisms of handlebodies.

We give precise definitions. A \emph{compression body} is a
manifold pair $(Q,F)$ obtained from a compact and orientable
surface $F$ in the following way. Consider the disjoint union of
the product $F\times I$ and finitely many balls (three dimensional
$0$--handles) $B$. Attach $1$--handles to
$(F\times\{1\})\cup\partial B$, obtaining $Q$. Identifying $F$
with $F\times\{0\}\subseteq Q$, we obtain the compression body
$(Q,F)$. Then $F\subseteq Q$ is the \emph{interior boundary} of
$(Q,F)$, denoted by $\partial_i Q$. The \emph{exterior boundary}
$\partial_e Q$ of $(Q,F)$ is the closure $\overline{\partial
Q-\partial_i Q}$. We allow empty or non-empty $\partial F$, but
$F$ cannot have sphere components. If $Q$ is homeomorphic to the
disjoint union of $F\times I$ with balls then $(Q,F)$ is said to
be \emph{trivial}.

A \emph{handlebody} $H$ is a connected compression body whose
interior boundary is empty, ie $H$ is obtained from attaching
$1$--handles to balls. The \emph{genus} of $H$ is the rank of
$\pi_1(H)$.

The following definition is due to Oertel:

\begin{definition}\label{D:red,irred}
An automorphism $f\co H\to H$ of a handlebody $H$ is said \emph{reducible}
if any of the following hold:
\begin{itemize}
   \item There exists an $f$--invariant (up to isotopy) nontrivial
   compression body $(Q,F)$ with $Q\subseteq H$, $\partial_e
   Q\subseteq \partial H$ and $F=\partial_i Q\neq\emptyset$ not
   containing $\partial$--parallel disc components.

   \item There exists an $f$--invariant (up to isotopy) collection of
   pairwise disjoint, incompressible, non-$\partial$--parallel and
   properly embedded annuli.

   \item $H$ admits an $f$--invariant (up to isotopy) $I$--bundle
   structure.
\end{itemize}

The automorphism $f$ is called \emph{irreducible} (or \emph{generic},
as in \cite{UO:Autos}) if
\begin{enumerate}
    \item $\partial f=f|_{\partial H}$ is pseudo-Anosov, and

    \item\label{I:closed reducing surface}
    there exists no \emph{closed reducing surface} $F$:
    a closed reducing surface is a surface $F\neq\emptyset$ which
    is the interior boundary $\partial_i Q$ of a nontrivial
    compression body $(Q,F)$ such that $Q\subseteq H$, $(Q,F)$
    is $f$--invariant (up to isotopy) and $\partial_e Q=\partial H$.
\end{enumerate}
\end{definition}

An obvious remark is that this definition of irreducible
automorphism excludes the periodic case. Also, a ``closed reducing
surface'' in \eqref{I:closed reducing surface} above is indeed
closed because the exterior boundary $\partial_e Q=\partial H$ is
closed (recall that the boundaries of $\partial_i Q$ and
$\partial_e Q$ coincide).

\begin{theorem}[Oertel \cite{UO:Autos}]\label{T:oertel classification}
An automorphism of a handlebody is either
\begin{enumerate}
    \item periodic,

    \item reducible, or

    \item irreducible.
\end{enumerate}
\end{theorem}

We note that the theorem above is not entirely obvious. For
example, one must show that if an automorphism $f\co H\to H$ of a
handlebody does not restrict to a pseudo-Anosov $\partial f$ on
$\partial H$, then $f$ is actually reducible according to
\fullref{D:red,irred}, or periodic.

Our interest is precisely in the irreducible case, which is in
many ways analogous to the pseudo-Anosov case for surfaces. An
important similarity is related to the existence of certain
invariant projective measured laminations (see \cite{UO:Autos}, and
\fullref{T:generic} and \fullref{Rs:laminations} below). A
good part of the original article is dedicated to the construction
of these laminations, which depends on many choices. Among these
one has to choose a handlebody $H_0\subseteq\mathring{H}$
``concentric'' with $H$, in the sense that the complement
$H-\mathring{H}_0$ is a product. Also, the automorphism $f\co H\to
H$ must be isotoped so that $H_1=f(H_0)$ contains $H_0$ in its
interior and $\bigcup_{i\in\mathbb{Z}} f^i(H_0)=\mathring{H}$.
This yields a nested chain
$H_i\subseteq\mathring{H}_{i+1}\subseteq H_{i+1}$, where
$H_i=f^i(H_0)$, $i\in\mathbb{Z}$. Also,
$\mathring{H}=\bigcup_{i\in\mathbb{Z}}H_i$. Next, a handle
decomposition $\cal{H}_0$ of $H_0$ as union of $0$ and
$1$--handles is needed (alternatively, one can choose a complete
system of discs, as in \cite{UO:Autos}). The 1--handles come with
a product structure $D^2\times I$.

We refer the reader to \cite{UO:Autos} for details on the
construction of the laminations. Important properties are
summarized in the following theorem and remarks. See
\fullref{Rs:laminations} for a comment on the singularities of the
one-dimensional lamination.

\begin{theorem}[Oertel]\label{T:generic}
Let $f\co H\to H$ be a generic automorphism of a handlebody. Then
there exist a two--dimensional measured lamination $(\Lambda,\mu)$
of $\ring{H}$ with full support, a one--dimensional measured
\emph{singular} lamination $(\Omega,\nu)$ in $H_0$ which is transverse
to $\Lambda$ and also with full support, $\widehat f$ isotopic to $f$
and $\lambda>1$ such that
\begin{enumerate}
    \item $\widehat f(\Lambda,\mu)=(\Lambda,\lambda\mu)$,

    \item $\widehat f(\Omega,\nu)=(\Omega,\lambda^{-1}\nu)$,

    \item\label{I:index} the leaves of $\Lambda$ are open discs and
    \emph{fill} $H_0$,
    in the sense that $\Lambda\cap H_0$ consists of
    essential discs in $H_0$ whose complement $H_0-\Lambda$
    consists of contractible components,

    \item $\Lambda\cup\partial H$ is closed in $H$,

    \item $\Lambda\cap\Omega$ is disjoint from the singular set of
    $\Omega$.
\end{enumerate}
\end{theorem}
\begin{obss}\label{Rs:laminations}
The handle decomposition $\cal{H}_0$ of $H_0$ determines handle
decompositions $\cal{H}_i$ of any $H_i=f^i(H_0)$,
$i\in\mathbb{Z}$, through $f^i$. Each 1--handle has a fixed
product structure $D^2\times I$. Consider the corresponding
product foliation by discs. We call a leaf of this foliation \emph{a disc
dual to the 1--handle}, or just \emph{a dual disc}. The
representative $\widehat f$ can be chosen so that the 1--handles of
$\cal{H}_1$ and $\cal{H}_0$ are \emph{compatible} in the sense that
1) dual discs in $H_1$ intersect $H_0$ in dual discs and 2)
$I$--fibers of the dual foliation of a 1--handle of $\cal{H}_0$ by
intervals intersect 1--handles of $\cal{H}_1$ in $I$--fibers. In
fact similar properties hold for any $\cal{H}_i$, $\cal{H}_j$,
$i<j$. For instance, dual discs in $H_j$ intersect $H_i$ in dual
discs.

Consider the intersection $\Lambda\cap H_0$ referred to in
property \eqref{I:index} in the theorem. It consists of a union of
discs dual to the 1--handles of the decomposition $\cal{H}_0$ (a
Cantor set of such discs). A similar description also holds for
any $H_i$ by the invariance of $\Lambda$ under $f$. For instance,
$\Lambda\cap H_1$ consists of families of discs essential in $H_1$
which laminate the 1--handles of $\cal{H}_1$ with dual discs.

Under the same point of view as above one can regard $\Lambda$ as
being obtained from $\Lambda_0=\Lambda\cap H_0$ by considering the
union of discs $\Lambda=\bigcup_{i\geq 0} f^i(\Lambda_0)$, where
$f^i(\Lambda_0)\subseteq f^j(\Lambda_0)$ whenever $i\leq j$.

The lamination $\Omega$ contains a \emph{singular set}
$S(\Omega)=\bigcap_{i\in\mathbb{Z}} f^i(h^0)$, where $h^0$ is the
union of $0$--handles of the decomposition $\cal{H}_0$. One can
choose $\widehat f$ so that $S(\Omega)$ is finite. The complement
$\Omega-S(\Omega)$ is a lamination of $H_0-S(\Omega)$. That is
where the measure $\nu$ is defined and is its support. Also, the
notion of tangency to $\Omega$ is not defined at $S(\Omega)$.
Therefore, in general, by saying that a surface $F$ is \emph{transverse}
to $\Omega$ we assume, in particular, that $F\cap
S(\Omega)=\emptyset$. The intersection of a 1--handle with
$\Omega$ consists of $I$--fibers of the 1--handle.
\end{obss}

The problem with \fullref{T:generic} is that the laminations,
whose construction depends on many choices, are not unique in any
reasonable sense. For example, the scalar $\lambda$, called \emph{the
growth rate of $f$ with respect to $\cal{H}$} or just \emph{the
growth rate} of $f$, which is a measure of the complexity of the
automorphism, is not unique. This phenomenon is not unlike what
happens with automorphisms of surfaces, if one allows the
invariant laminations to have monogons for complementary
components.

The main problem we shall address here, though not solve, is the
following.

\begin{problem*}
Characterize canonical invariant laminations for a given
irreducible automorphism of a handlebody.
\end{problem*}

In the case of surfaces, a solution is to fix a hyperbolic metric
in its interior and work with certain geodesic laminations which,
among other properties, realize minimal growth. For automorphisms
of handlebodies there is no such solution. Still, there is a
minimum in the set of possible growth rates. Naturally, canonical
laminations must yield minimal growth, ie the corresponding
automorphism must be most \emph{efficient} in its isotopy class. A
step (a big step, we believe) in the direction of solving the
problem above would then be to characterize minimal growth. Oertel
gives a necessary condition \cite{UO:Autos}.

\begin{theorem}[Oertel]\label{T:incompressibility}
If $\lambda$ is minimal then $\Lambda$ has the \emph{incompressibility
property}: for each leaf $L$ of $\Lambda$ the
complement $L-\mathring{H}_0$ is incompressible in $\ring
H-\mathring{H}_0$.
\end{theorem}

A clear sufficient condition, much stronger than
incompressibility, is that the leaves of $\Omega$ do not
``back-track'' (with respect to $\Lambda$) in $H_0$ (see
the work of Bestvina and Handel \cite{BH:Tracks}). Not surprisingly it
cannot always be realized
(see remark below).

We will also consider the problem of constructing examples of irreducible
automorphisms of handlebodies. As in any field of mathematics,
examples provide a useful investigative tool. The construction of
irreducible automorphisms of handlebodies is not an obvious task. The
main difficulty resides in proving that a given example does not admit
closed reducing surfaces, see \fullref{D:red,irred} (the other property,
that the restriction to the boundary is pseudo-Anosov, can be achieved
with the help of some well-established tools; see the articles by Penner
\cite{RP:88} and by Bestvina and Handel \cite{BH:Surfaces}).

A result of Bonahon \cite{FB:RibbonKnots} implies that any automorphism
of a genus two handlebody whose restriction to the boundary is
a pseudo-Anosov automorphism is then irreducible (see also the
article by Long \cite{DL:Bounding} and the author's doctoral thesis
\cite{LC:thesis}). This is not true for higher genus handlebodies.

\begin{obs}\label{R:examples}
Bonahon's result may be used to generate interesting examples. For
instance, in \cite{AFFL:Heegaard}, Fathi and Laudenbach build an
automorphism of a genus two handlebody which (1) restricts to the
boundary as a pseudo-Anosov automorphism --- thus, as mentioned
before, is irreducible --- and (2) induces the identity on the
fundamental group. In particular, the leaves of $\Omega$ have to
``back-track''. Such an example illustrates the richness of
irreducible automorphisms of handlebodies when compared with
pseudo-Anosov automorphisms of surfaces, whose complexity is
captured on the level of the fundamental group.
\end{obs}

\subsection{Summary of results}

Our main results address the problem of characterizing minimal
growth of a given irreducible automorphism of a handlebody. We
will identify a property on some two--dimensional laminations
$\Lambda$, which we call ``tightness''\footnote{In fact, being
tight is a property of the pair of measured laminations.} in
\fullref{D:tight}. Essentially, $\Lambda$ is ``tight'' if
the weighted intersection of its leaves with $(\Omega,\nu)$ is
minimal. The property of being tight is (strictly) stronger than
that of Oertel's incompressibility and (strictly) weaker than that
of having ``no back-tracking''. We conjecture that it
characterizes minimal growth.

\begin{conj-equiv}
The growth rate $\lambda$ is minimal if and only if $\Lambda$ is
\tight.
\end{conj-equiv}

This work will prove one direction:

\begin{thm-minim}
If $\Lambda$ is tight then $\lambda$ is minimal.
\end{thm-minim}

As for necessity, the problem is harder. We will prove it only
under some technical hypotheses (\fullref{P:valence}). These
hypotheses are useful: we will show that they can be assumed for
genus two handlebodies. In this case tightness characterizes
minimal growth:

\begin{corr-true}
\fullref{Cj:equivalence} is true for genus two
handlebodies.
\end{corr-true}

Moreover, tightness yields results concerning the growth rates. We prove:

\begin{corr-small}
If $\Lambda$ is  tight, then the growth rate $\lambda$ (which is
minimal) is less than or equal to the growth rate of the
restriction of the automorphism $f\co H\to H$ to the boundary
$\partial H$ (which is pseudo-Anosov).
\end{corr-small}

The following is a corollary of \fullref{T:minimality}.

\begin{corr-growth}
If $\Lambda$ is \tight\ then the minimal growth
$\lambda_{min}(f^n)$ of any power $f^n$ is
$\left(\lambda_{min}(f)\right)^n$.
\end{corr-growth}

From our point of view the measures $\nu$, $\mu$ on the invariant
laminations $\Omega$, $\Lambda$ and the corresponding growth rate
$\lambda$ come from eigenvectors and eigenvalue of certain
incidence matrices of a handle decomposition $\cal{H}$ of the
handlebody (or complete disc system $\cal{E}$) obtained through
$f$. The construction depends on such a matrix $M$ being
non-negative and \emph{irreducible}. That means that for any $1\leq
i,j\leq \text{dim}(M)$ there exists a power $M^n$, $n\geq 1$, in
which the corresponding entry $ij$ is not zero. If the
automorphism $f$ is irreducible we can assume that the matrix $M$
is irreducible (see the paper \cite{UO:Autos} by Oertel), which is
required in the construction of the laminations. Irreducible matrices have
many nice spectral properties. For instance its spectral radius, which
we denote by $\lambda(M)$, is realized by a positive real eigenvalue. We
call it the \emph{Perron--Frobenius eigenvalue of $M$}. The corresponding
positive eigenvector, which is well-defined (up to scaling, naturally),
is called the \emph{Perron--Frobenius eigenvector}. The following result
is useful (see, for example, the articles \cite{ES:73} by Seneta, and
\cite{BH:Tracks} by Bestvina and Handel).

\begin{proposition}\label{P:subinvariance}
Let $M$ be a non-negative and irreducible $n\times n$ matrix and
$v\in \mathbb{R}^n$ with $v_i\geq 0$ and $v\neq 0$. If
$$(Mv)_i\leq \lambda v_i \quad \text{for all } i$$
then $\lambda(M)\leq\lambda$ and $v_i>0$. If, moreover,
$(Mx)_i<\lambda x_i$ for some $i$, then $\lambda(M)<\lambda$.
\end{proposition}

We will also present a method for generating a certain class of
irreducible automorphisms. This method produces examples on higher
genus handlebodies. The fact that our techniques fail to prove
\fullref{Cj:equivalence} for a general automorphism of a handlebody of
genus greater than two makes it especially important to study examples
in higher genus cases. Our method will depend on the two following
results. See \fullref{D:penner pair} and \fullref{R:irreducible} below
for important definitions or the work of Penner \cite{RP:88} and of
Bestvina and Handel \cite{BH:Tracks}, respectively, for more details.

\begin{theorem}[Penner]\label{T:penner}
Let $\mathcal{C}$, $\mathcal{D}$ be two systems of closed curves
in an oriented surface $S$ with $\chi(S)<0$. Assume that $\cal{C}$
and $\cal{D}$ intersect efficiently, do not have parallel
components and fill $S$.  Let $f\co S\to S$ be a composition of
Dehn twists: right twists along curves of $\mathcal{C}$ and left
twist along curves of $\mathcal{D}$. If a twist along each curve
appears at least once in the composition, then $f$ is isotopic to
a pseudo-Anosov automorphism of $S$.
\end{theorem}

\begin{theorem}\label{T:pA irreducible}
Let $S$ be a compact surface with $\chi(S)<0$ and precisely one
boundary component. An automorphism $f\co S\to S$ is
pseudo-Anosov if and only if the map
$f_*^n\co\pi_1(S)\to\pi_1(S)$ is
\emph{irreducible} for all $n>0$.
\end{theorem}

\begin{obs}\label{R:irreducible}
We recall from Bestvina--Handel \cite{BH:Tracks} the definition
of an irreducible automorphism of a free group $F$. The \emph{outer
automorphism group} of $F$ is obtained from the group of automorphisms
(ie self-isomorphisms) of $F$ by identifying any two which differ by an
inner isomorphism. An ``outer automorphism'' (ie an equivalence class)
$[\varphi]$ is said \emph{reducible} if the following holds. There are
proper free factors $F_1,\dots,F_k$ of $F$ such that $[\varphi]$ permutes
the conjugacy classes of the $F_i$'s and such that $F_1\ast\cdots\ast
F_k$ is a free factor (not necessarily proper) of $F$. If $[\varphi]$
is not reducible it is said \emph{irreducible}.

The following abuse is present in this paper. When considering an
automorphism of a group (typically of a fundamental group) it will
often be regarded as its outer automorphism class. For instance,
$f_*^n$ makes sense in the statement of \fullref{T:pA
irreducible} as an outer automorphism. As such it makes sense to
wonder whether it is reducible or irreducible (note that
$\pi_1(S)$ is free).
\end{obs}

The results of this article are divided in two following sections.
In \fullref{S:examples} we will describe our method for
generating examples of irreducible automorphisms. We will develop
a particular case and then generalize it in \fullref{T:method}. Its
statement depends on some technical
constructions unsuited for this introduction. We will then use it
to build a certain irreducible automorphism of a genus four
handlebody (\fullref{E:method}). It will help in motivating
the relevance of the property of tightness. For this reason we
shall determine a certain pair of invariant laminations for this
automorphism and estimate the corresponding growth rate. The
two--dimensional lamination will have Oertel's incompressibility
property.

\fullref{S:tightness} is dedicated to the tightness property.
We shall see that the example built in the preceding section does
not realize minimal growth. The lack of tightness, which we will
define then, will appear naturally there. This will be done
through the existence of ``tightening discs'', which will be our
main objects in dealing with lack of tightness. In the remainder
of the section we shall prove the theorems and corollaries on
tightness already mentioned.

We adopt the following notations and conventions. Given a topological
space $A$ (typically a manifold or sub-manifold), $\overA$ denotes
its topological closure, $\mathring{A}$ its interior and $|A|$ its
number of connected components. If $H$ is a handlebody we denote a
handle decomposition of $H$ by $\cal{H}$.  By considering co-cores
of $1$--handles (which we may also call \emph{dual discs}, see
\fullref{Rs:laminations}) it is clear that a handle decomposition
$\cal{H}$ of $H$ corresponds to a complete system of discs
$\cal{E}\subseteq H$. In fact, the set of handle decompositions
and the set of complete disc systems are, up to isotopy, in $1{-}1$
correspondence. This remark is relevant for while the paper by Oertel
\cite{UO:Autos} uses discs systems --- because the author focuses on the
two--dimensional lamination --- we shall use handle decompositions ---
because we focus more on the one--dimensional lamination.  The incidence
matrix associated to a $\cal{H}$ corresponds to the transpose of the
incidence matrix associated to the corresponding $\cal{E}$. There is also
an embedded graph $\Gamma\subseteq H$ dual to $\cal{E}$, with vertices
corresponding to $0$--handles of $\cal{H}$ and edges corresponding
to $1$--handles. We can then regard $H$ as a fibered neighborhood of
$\Gamma$. Such embedded graphs are also, up to isotopy, in $1{-}1$
correspondence with handle decompositions and complete disc systems.

I thank Ulrich Oertel for his helpful suggestions in his role as
dissertation advisor, and also for laying the foundations on which
the research in this paper is built.  I thank the referee for
reading the paper carefully, and for suggesting improvements to
the exposition.

This research was partly supported by CNPq---Brazil,
CAPES---Brazil fellowship BEX0292/99-0 and FAPESP---Brazil
fellowship 03/06914-2.

%%%%%%%%%%%%%%%%%%%%%%%%%%%%%%%%%%%%%%%%%%%%%%%%%%%%%%%%%%%%%%%%%%%%%%%%%%%%%%%%%%%%%%%%%%%%%%%%%%%%%%%%%
%%%%%%%%%%%%%%%%%%%%%%%%%%%%%%%%%%%%%%%%%%%%%%%%%%%%%%%%%%%%%%%%%%%%%%%%%%%%%%%%%%%%%%%%%%%%%%%%%%%%%%%%%
%%%%%%%%%%%%%%%%%%%%%%%%%%%%%%%%%%%%%%%%%%%%%%%%%%%%%%%%%%%%%%%%%%%%%%%%%%%%%%%%%%%%%%%%%%%%%%%%%%%%%%%%%
%%%%%%%%%%%%%%%%%%%%%%%%%%%%%%%%%%%%%%%%%%%%%%%%%%%%%%%%%%%%%%%%%%%%%%%%%%%%%%%%%%%%%%%%%%%%%%%%%%%%%%%%%

\section{Examples}\label{S:examples}

\subsection{An example}\label{SS:first example}

We will develop a simple particular case of the method that will
be obtained in the following subsection. Let $H$ be a genus 2
handlebody. An automorphism of $H$ will be described as a
composition of Dehn twists along two annuli and a disc. We shall
prove that it is irreducible by showing that its restriction to
$\partial H$ is pseudo-Anosov and that, for an algebraic reason,
there can be no closed reducing surface. The argument that proves
this last part is distinct from that of Bonahon \cite{FB:RibbonKnots} for
genus 2 handlebodies and, with the right hypotheses, generalizes
to higher (even) genus handlebodies.

%%%%%%%%%%%%%%%%%%%%%%%%%%%%%%%%%%%%%%%%%%%%%%%%%%
%
%   Exemplo 1
%
%%%%%%%%%%%%%%%%%%%%%%%%%%%%%%%%%%%%%%%%%%%%%%%%%%

\begin{example}\label{E:first example}
We first construct a pseudo-Anosov automorphism $\varphi\co S\to
S$ of the once punctured oriented torus $S$. It will be defined as
a composition of Dehn twists along two curves.

We represent $S$ as a cross with pairs of opposite sides
identified as shown in \fullref{F:toro}. Fixing a base point in
$S$ we note that $\pi_1(S)$ is the free group on two generators.

\begin{figure}[ht!]
\labellist\small
\pinlabel {$\alpha_0$} [b] at 356 367
\pinlabel {$\alpha_1$} [tl] at 330 76
\endlabellist
\centerline{\includegraphics[width=5cm]{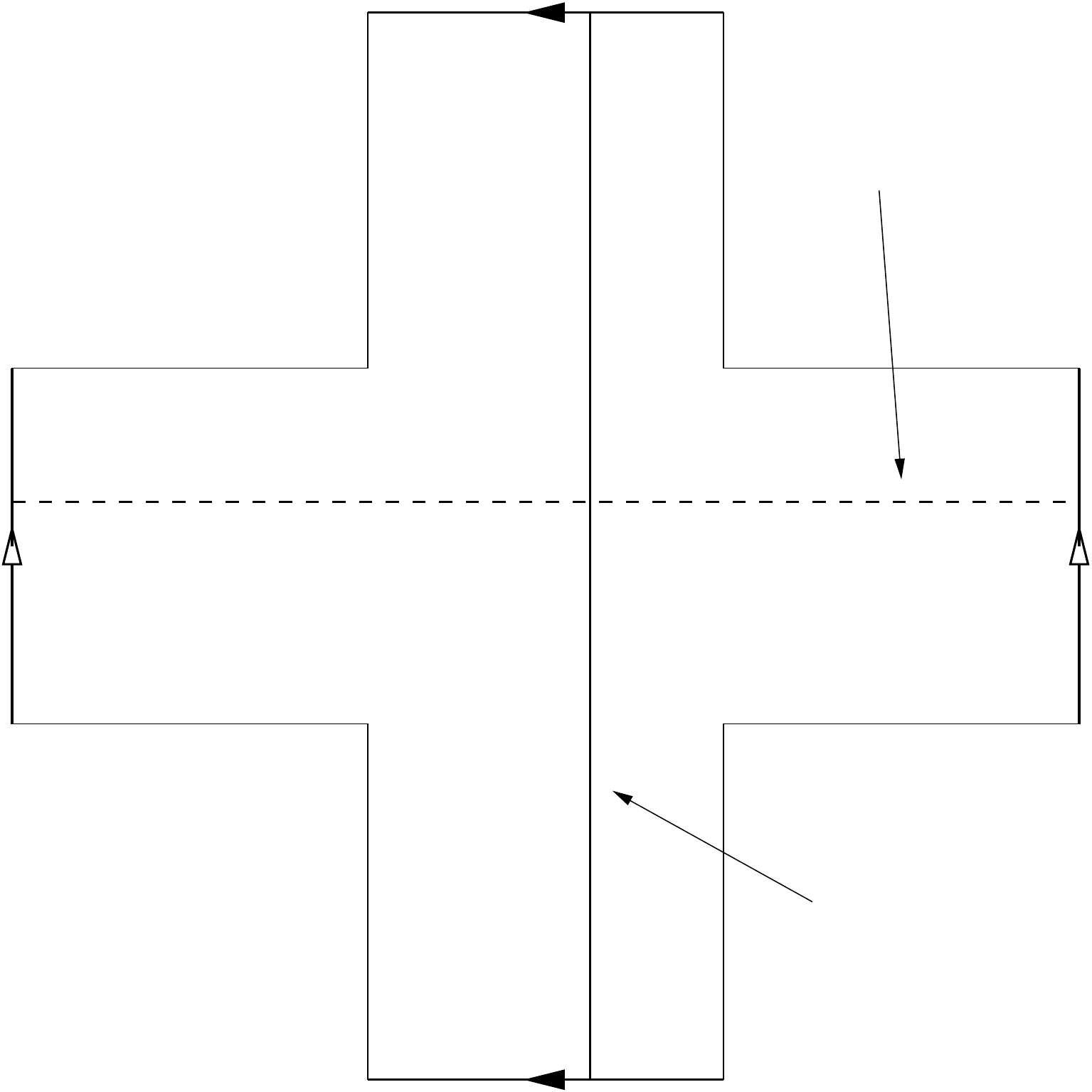}}
\caption{The oriented surface $S$ and the curves
$\alpha_0$, $\alpha_1$}
\label{F:toro}
\end{figure}

Let $\alpha_0$, $\alpha_1$ be simple closed curves as in the
figure. It is easy to verify that the systems $\cal
C=\{\alpha_0\}$ e $\cal D=\{\alpha_1\}$ satisfy the hypotheses of
\fullref{T:penner} (Penner). Let $T_0^-$ be the left Dehn
twist along $\alpha_0$ and $T_1^+$ the right twist along
$\alpha_1$. We define
$$\varphi=T_1^+\circ T_0^-.$$
By \fullref{T:penner} $\varphi$ is pseudo-Anosov. And by
\fullref{T:pA irreducible} any positive power $\varphi_*^n$ of
the induced isomorphism $\varphi_*\co\pi_1(S)\to\pi_1(S)$ is
\emph{irreducible} (see \fullref{R:irreducible}). We note this fact
for future use.

We now consider the handlebody $H=S\times I$, identifying $S$ with
$S\times\{1\}$. The orientation of $S$ then determines an
orientation on $H$ through inclusion. Now lift $\varphi$ to $H$,
obtaining $\phi\co H\to H$, a composition of twists along the
annuli $A_0=\alpha_0\times I$, $A_1=\alpha_1\times I$ as in
\fullref{F:aneis}. Identifying $\pi_1(H)$ with $\pi_1(S)$
yields $\phi_*=\varphi_*$.

\begin{figure}[ht!]
\labellist\small\hair1.5pt
\pinlabel {$\Delta$} [b] at 960 382
\pinlabel {$\mathcal{C}$} [b] at 699 384
\pinlabel {$\mathcal{D}$} [tl] at 963 66
\pinlabel {$A_0$} [bl] at 395 345
\pinlabel {$A_1$} [l] at 371 75
\pinlabel {$H$} at 71 417
\endlabellist
\centerline{\includegraphics[scale=0.30]{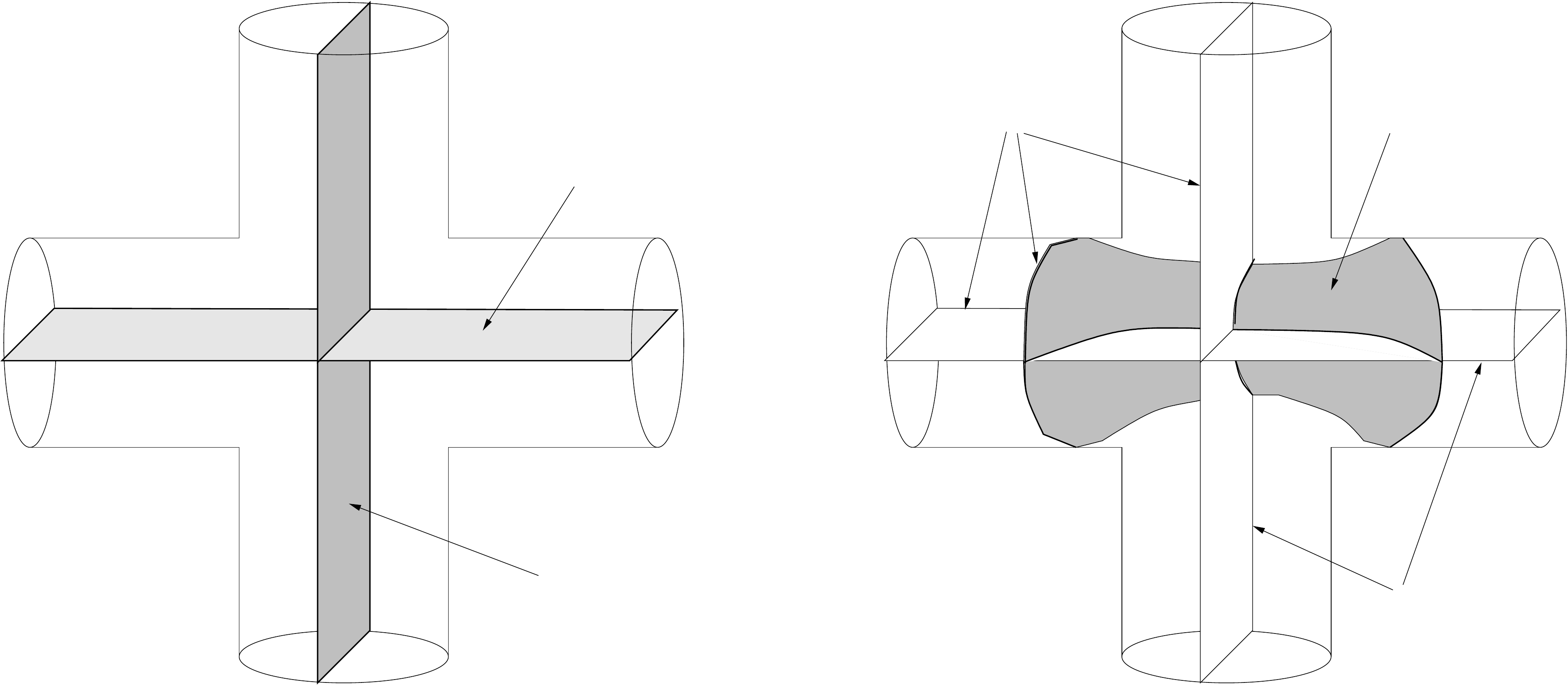}}
\caption{The automorphism $f$ is defined as a composition of Dehn twists along
the annuli $A_0$, $A_1$ and the disc $\Delta$.}
\label{F:aneis}
\end{figure}
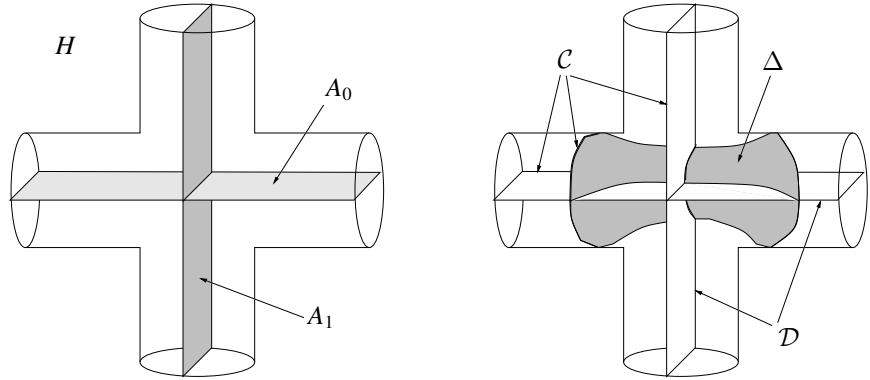

Finally, we will obtain the desired irreducible automorphism $f\co
H\to H$ by composing $\phi$ with a twist along a disc $\Delta$,
shown in \fullref{F:aneis}. Let $T_\Delta^+$ be the right Dehn
twist along $\Delta$. We define:
$$
f=T_\Delta^+\circ\phi.
$$

\begin{proposition}\label{P:reducao bordo}
The automorphism $f\co H\to H$ is irreducible.
\end{proposition}

The lack of closed reducing surfaces will come from the following
general lemma. See \fullref{R:irreducible} for the definition
of a \emph{reducible automorphism of a free group}.

\begin{lemma}\label{L:reducible pi_1}
Let $g\co H\to H$ be an automorphism of a handlebody $H$ such that
$\partial g$ is pseudo-Anosov. If $g$ is reducible then, for some
$n\in\mathbb{N}$, $g_*^n\co\pi_1(H)\to\pi_1(H)$ is reducible.
\end{lemma}
\begin{obs}\label{R:base point}
It is clear that the result does not depend on the choice of base
point.
\end{obs}
\begin{proof}
Since $\partial g$ is pseudo-Anosov and $g$ is reducible it
follows from \fullref{D:red,irred} and \fullref{T:oertel classification}
that there exists a $g$--invariant
and nontrivial compression body $Q$ such that $\partial_e
Q=\partial H$ and $\partial_i Q\neq\emptyset$. Let
$F\subseteq\partial_iQ$ be a component of the closed reducing
surface $\partial_i Q$. By removing 1--handles of $Q$ it is easy
to see that $F$ bounds a handlebody $J\subseteq\ring H$. By
choosing a base point in $J$ (see the remark above) and omitting
the obvious inclusion homomorphisms we claim that
$$
\pi_1(H)=\pi_1(J)*G,
$$
where $G$ is not trivial. To see this first consider the connected
and nontrivial compression body $Q'=\overline{H-J}$, whose
boundary decomposes as $\partial_i Q'=F$ and $\partial_e
Q'=\partial H$. A compression body structure of $Q'$ gives it as a
product $F\times I$ to which 1--handles are attached. Regarding
$F\times I\subseteq Q'\subseteq H$ we see that the handlebody
$J'=(F\times I)\cup J$ deformation retracts to $J$ (so
$\pi_1(J')=\pi_1(J)$ through inclusion). Using the compression
body structure of $Q'$ we can regard $H$ as obtained from $J'$ by
attaching 1--handles to $\partial J'$. Since $\partial J'$ is
connected, we can moreover assume that these 1--handles are
attached to a disc in $\partial J'$, which gives
$\pi_1(H)=\pi_1(J')*G=\pi_1(J)*G$, where $G$ is a free group
(whose rank equals the minimal number of 1--handles in a
compression body structure of $Q'$). But $Q'$ is not trivial, so
$G$ is not trivial, proving the claim. Therefore $\pi_1(J)$ is a
proper free factor of $\pi_1(H)$.

Let $g^n$ be the first power of $g$ preserving $J$. Isotoping $g$
we assume moreover that the base point is fixed by $g^n$. From
$$
g^n(J)=J
$$
it follows that $g_*^n(\pi_1(J))$ is conjugate to $\pi_1(J)$,
hence $g_*^n$ is reducible.
\end{proof}

\begin{proof}[Proof of \fullref{P:reducao bordo}]
We need to prove that $\partial f=f|_{\partial H}$ is
pseudo-Anosov and that $f$ does not admit closed reducing
surfaces.

We start by verifying that $\partial f$ is pseudo-Anosov. It is
given as a composition of Dehn twists: left twists along curves of
$$
\mathcal{C}=\left\{\,(\alpha_0\times\{1\})\,,\,(\alpha_1\times\{0\})\,\right\},
$$
(see \fullref{F:aneis}) and right twists along curves of
$$
\mathcal{D}=\left\{\,(\alpha_0\times\{0\})\,,\,(\alpha_1\times\{1\})\,,\,\partial\Delta\,\right\}.
$$
We now note that $\mathcal{C}$, $\mathcal{D}$ satisfy the
hypotheses of \fullref{T:penner} (also see \fullref{D:penner pair}),
hence $\partial f$ is pseudo-Anosov.

Now suppose by contradiction that there exists a closed reducing
surface. By \fullref{L:reducible pi_1} there exists $n$ such
that $f_*^n$ is reducible. But $f=(T_\Delta^+)\circ\phi$ and the
twist $T_\Delta^+$ (along a disc) induces the identity in
$\pi_1(H)$. Therefore, recalling that $\pi_1(H)$ is identified
with $\pi_1(S)$, we have that $f_*^n=\phi_*^n=\varphi_*^n$, which
was seen before to be irreducible for any $n$, a contradiction.
This shows that there are no closed reducing surfaces, completing
the proof.
\end{proof}
\end{example}

%%%%%%%%%%%%%%%%%%%%%%%%%%%%%%%%%%%%%%%%%%%%%%%%%%%%%%%%%%%%%%%%%%%%%%%%%%%%%%%%%%%%%%%%%%
%%%%%%%%%%%%%%%%%%%%%%%%%%%%%%%%%%%%%%%%%%%%%%%%%%%%%%%%%%%%%%%%%%%%%%%%%%%%%%%%%%%%%%%%%%

\subsection{A method for generating irreducible
automorphisms} \label{SS:method}

The construction of \fullref{E:first example} may be
generalized to provide a method for generating a larger class of
irreducible automorphisms of handlebodies (Theorems \ref{T:method}
and \ref{T:method2}). It partially solves a problem proposed by Oertel
\cite{UO:Autos}.

\begin{definition}\label{D:penner pair}
Let $(\cal{C},\cal{D})$ be a pair of curve systems in a compact,
connected and orientable surface $S$ with $\chi(S)<0$. It is
called a \emph{Penner pair in $S$} if $\cal C$, $\cal D$ satisfy
the hypotheses of Penner's \fullref{T:penner}, that is,
\begin{enumerate}
    \item each $\cal C$, $\cal D$ is a finite collection of simple,
    closed and pairwise disjoint essential curves without parallel
    copies,

    \item $\cal{C}$ and $\cal{D}$ intersect efficiently, do not
    have parallel components and \emph{fill} $S$
    (ie the components of $S-(\cal{C}\cup\cal{D})$ are either
    contractible or deformation retract to a component of $\partial S$).
\end{enumerate}
Suppose that $(\cal C,\cal D)$ is a Penner pair. An automorphism
$\varphi$ of $S$ obtained from $\cal C$, $\cal D$ as in
\fullref{T:penner} is called a \emph{Penner automorphism
subordinate to $(\cal C,\cal D)$} (which is in particular
pseudo-Anosov).

If $\partial S\neq\emptyset$ then a properly embedded and
essential arc $\theta$ is called \emph{dual to $(\cal C,\cal D)$}
if $\theta$ intersects $\cal C\cup\cal D$ transversely and in
exactly one point $p\notin\cal C\cap\cal D$.
\end{definition}

We constructed the irreducible automorphism in \fullref{E:first example}
by lifting a pseudo-Anosov automorphism of
a surface to a product and composing it with a twist on a disc.
The general method will be similar. Our interest in dual arcs is
that we can use them to construct discs that will yield
irreducible automorphisms.

Throughout this subsection we fix a compact, connected and
oriented surface $S$ with $\partial S\neq\emptyset$ and define
$H=S\times I$, which is a handlebody. We identify $S$ with
$S\times\{1\}\subseteq H$, inducing orientation in $H$. We also
fix a base point in $S\times\{1\}$ for both $S$ and $H$ and
identify $\pi_1(H)$ with $\pi_1(S)$.

Given a Penner pair $(\cal C,\cal D)$ in $S$ and a dual arc
$\theta$ we build a disc $\Delta_\theta$ in $H$ in the following
way. Let $\gamma$ be the curve of $(\cal C,\cal D)$ that $\theta$
intersects and assume without loss of generality that
$\gamma\subseteq\cal{C}$. Let $D=\theta\times I\subseteq H$. Then
$\partial D$ intersects $\gamma_1=\gamma\times\{1\}$ in a point.
Now let $\Delta_\theta$ be the \emph{band sum} of $D$ with itself
along $\gamma_1$.  This means that $\Delta_\theta$ is obtained
from $D$ and $\gamma_1$ by the following construction. Consider a
regular neighborhood $N=N(D\cup\gamma_1)$. Then
$\Delta_\theta=\overline{\partial N-\partial H}$ is a properly
embedded disc.

\begin{theorem}\label{T:method}
Suppose that $\partial S\neq\emptyset$ has exactly one component.
Let $(\cal C,\cal D)$ be a Penner pair in $S$ with dual arc
$\theta$ and $\varphi\co S\to S$ a Penner automorphism subordinate
to $(\cal C,\cal D)$. Let $\widehat\varphi\co H\to H$ be the lift of
$\varphi$ to the product $H=S\times I$ and $\Delta_\theta\subseteq
H$ the disc constructed from the arc $\theta$ as above. Then there
exists a simple Dehn twist $T_{\Delta_\theta}\co H\to H$ along
$\Delta_\theta$ such that the composition
$$
\widehat\varphi\circ T_{\Delta_\theta} \co H\to H
$$
is an irreducible automorphism of $H$.
\end{theorem}

The key to the proof is the verification that $\cal{C}$, $\cal{D}$
and $\partial\Delta_\theta$ determine a Penner pair in $\partial
H$.

\begin{lemma}\label{L:twisting disc}
Let $S$, $(\mathcal{C},\mathcal{D})$, $\theta$, $H=S\times I$ and
$\Delta_\theta$ be as in the statement of \fullref{T:method}.
Let $\cal C_i=\cal C\times\{i\}\subseteq S_i=S\times\{i\}$ and
$\cal D_i=\cal D\times\{i\}\subseteq S_i=S\times\{i\}$, defining
$\cal C_0$, $\cal D_0\subseteq S_0$ and $\cal C_1$, $\cal
D_1\subseteq S_1$. Under these conditions the system
\begin{align}
\mathcal{Q} =& \cal D_0 \cup \cal C_1
\cup\{\,\partial\Delta_\theta\,\},\notag\\
\mathcal{R} =& \cal C_0 \cup \cal D_1 \notag,
\end{align}
of curves in $\partial H$,
determines a Penner pair $(\mathcal{Q},\mathcal{R})$ in $\partial
H$.
\end{lemma}
\begin{proof}
We start by making the obvious remarks that $\cal C_0$, $\cal
D_0$, $\cal C_1$, $\cal D_1\subseteq\partial H$ and
$\cal{C}_0\cap\cal{D}_1=\emptyset$,
$\cal{D}_0\cap\cal{C}_1=\emptyset$. Recall that we are assuming
that
$\theta\cap(\cal{C}\cup\cal{D})\subseteq\gamma\subseteq\cal{C}$.
We verify that
\begin{itemize}

    \item $\partial\Delta_\theta\cap\cal{D}_0=\emptyset,$ because
    $\left(\theta\times\{0\}\right)\cap\cal{D}_0=\emptyset$ and
    $\partial\Delta_\theta\cap S_0$ consists of two arcs in $S_0$
    parallel to $\theta\times\{0\}$,

    \item $\partial\Delta_\theta\cap\cal{C}_1=\emptyset,$ because
    $\partial\Delta_\theta\cap\gamma_1=\emptyset$ by construction.
\end{itemize}
Therefore each $\mathcal{Q} =\cal D_0 \cup \cal
C_1\cup\{\partial\Delta\}$ and $\mathcal{R} = \cal C_0 \cup\cal
D_1$ is a system of simple closed curves essential in $\partial
H$. To conclude that $(\mathcal{Q},\mathcal{R})$ is indeed a
Penner pair it remains to verify that $\mathcal{Q}\cup\mathcal{R}$
fills $\partial H$.

A component of $S-(\cal C\cup\cal D)$ either is a disc or an
annulus that retracts to $\partial S$. Therefore a component of
$\partial H-(\cal C_0\cup\cal D_0\cup \cal C_1\cup\cal D_1)$
either is a disc or an annulus $A$ (that retracts to $\partial
S\times I$). But $A\cap\partial\Delta_{\theta}$ is a union of four
arcs essential in $A$, hence each component of $\partial
H-(\mathcal{Q}\cup\mathcal{R})$ is a disc. In other words
$\mathcal{Q}\cup\mathcal{R}$ fills $\partial H$, completing the
proof.
\end{proof}

Instead of proving \fullref{T:method} we will prove the more
general result below, which clearly implies the other. We note
that twists on curves of $\cal C$, $\cal D$ in $S$ lift to twists
along annuli in $H$. We denote these systems of annuli by
$\widehat{\cal C}$, $\widehat{\cal D}$ respectively. The ``direction of a
twist'' along these vertical annuli should be understood as the
direction of its restriction to $S\times\{1\}\subseteq\partial H$.

\begin{theorem}\label{T:method2}
Let $(\mathcal{C},\mathcal{D})$, $S$, $\theta$, $H$ and
$\Delta_\theta$ be as in \fullref{T:method}. Let $f$ be a
composition $f\co H\to H$ of twists along the annuli of $\widehat{\cal
C}$, $\widehat{\cal D}$ and the disc $\Delta_{\theta}$: in one
direction along the annuli in $\widehat{\cal D}$ and in the opposite
direction along the annuli in $\widehat{\cal C}$ and the disc
$\Delta_{\theta}$. If each of these twists appear in the
composition at least once $f$ is irreducible.
\end{theorem}
\begin{proof}
We first show that $f_*^n\co\pi_1(H)\to\pi_1(H)$ is an irreducible
automorphism of a free group for any $n\geq 0$ (hence there can be
no closed reducing surface by \fullref{L:reducible pi_1}) and
then that $\partial f=f|_{\partial H}$ is pseudo-Anosov, thus
completing the proof that $f$ is irreducible.

Recall that $S$ is identified with $S\times\{1\}\subseteq H$ and
$\pi_1(S)$ with $\pi_1(H)$. Let $T_{\Delta_{\theta}}$ be a twist
along $\Delta_{\theta}$. Since
$(T_{\Delta_\theta})_*\co\pi_1(H)\to\pi_1(H)$ is the identity
($\Delta_\theta$ is a disc) the hypotheses on $f$ imply that
$f_*=\varphi_*$ for some Penner automorphism $\varphi\co S\to S$
subordinate to $(\cal C,\cal D)$. Penner automorphisms are
pseudo-Anosov so, given that $\partial S$ has a single component,
it follows from \fullref{T:pA irreducible} that $\varphi_*^n$
is an irreducible automorphism of $\pi_1(S)$ for any $n\geq 0$.
Therefore $f_*^n\co\pi_1(H)\to\pi_1(H)$ is irreducible, and then
$f$ does not admit closed reducing surfaces
(\fullref{L:reducible pi_1}).

To see that $\partial f$ is pseudo-Anosov, let
$(\mathcal{Q},\mathcal{R})$ be as in \fullref{L:twisting disc},
therefore a Penner pair. By construction the twists that compose
$f$ restrict to $\partial H$ as twists along curves of $\cal{Q}$
or $\cal{R}$. It is then straightforward to verify that $\partial
f$ is a Penner automorphism subordinate to $(\cal{Q},\cal{R})$,
hence pseudo-Anosov, completing the proof that $f$ is
irreducible.
\end{proof}

\begin{obs}\label{R:even genus}
Note that the conditions that $S$ is orientable and $|\partial
S|=1$ imply that $H$ has even genus.
\end{obs}

\begin{example}\label{E:method}
Consider $S$ a genus 2 surface minus a disc, represented in
\fullref{F:twisting curves} as an octagon whose sides are
identified according to the arrows.

\begin{figure}[ht!]
\labellist
\hair1pt
\pinlabel {$\alpha$} [l] at 450 398
\pinlabel {$\beta$} [r] at 27 152
\pinlabel {$\gamma$} [l] at 469 317
\pinlabel {$\delta$} [tl] at 452 44
\pinlabel {$\theta$} [l] at 937 209
\pinlabel {$A$} [b] at 257 403
\pinlabel {$A$} [r] at 64 214
\pinlabel {$B$} [br] at 120 350
\pinlabel {$B$} [tr] at 118 87
\pinlabel {$C$} [tl] at 382 84
\pinlabel {$C$} [bl] at 382 350
\pinlabel {$D$} [t] at 255 27
\pinlabel {$D$} [l] at 440 214
\endlabellist
\centerline{\includegraphics[scale=0.35]{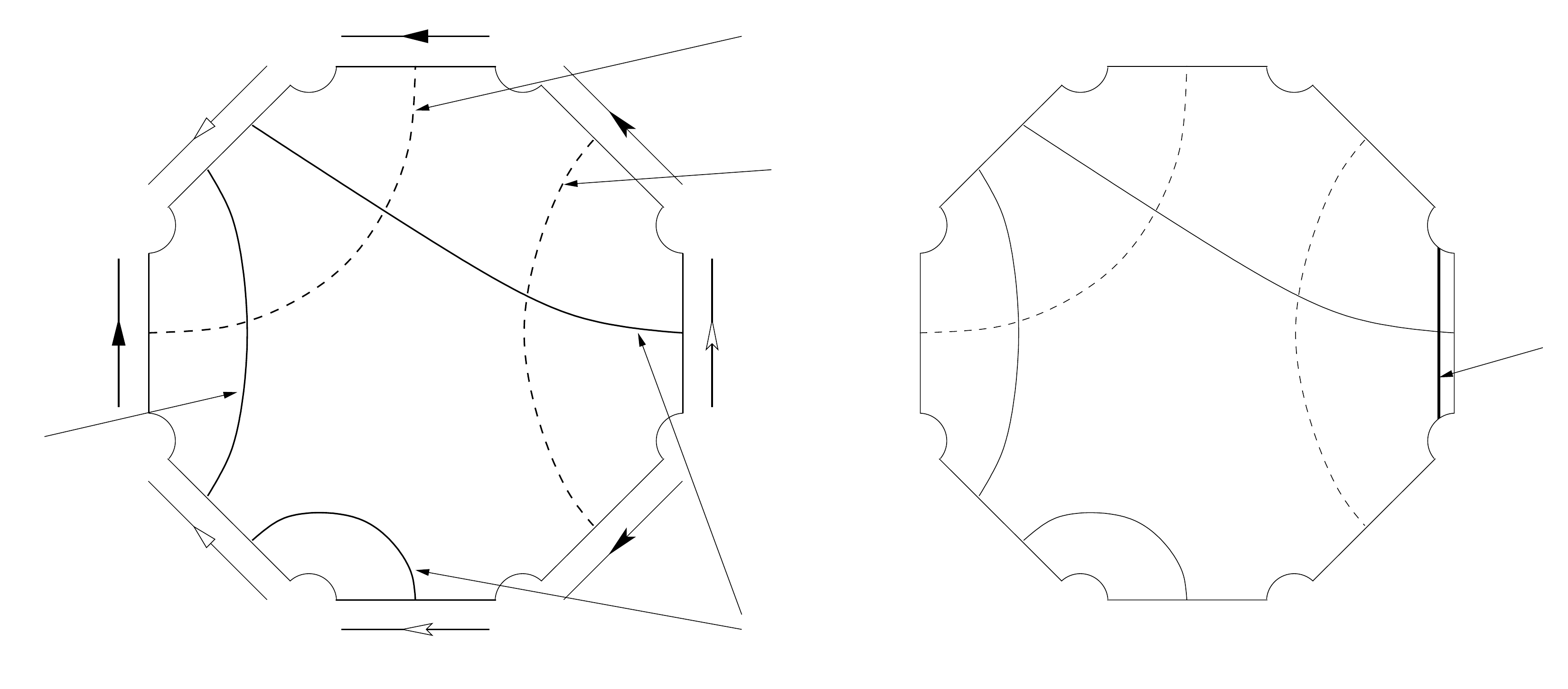}}
\caption{A Penner pair in $S$, with dual arc $\theta$}
\label{F:twisting curves}
\end{figure}
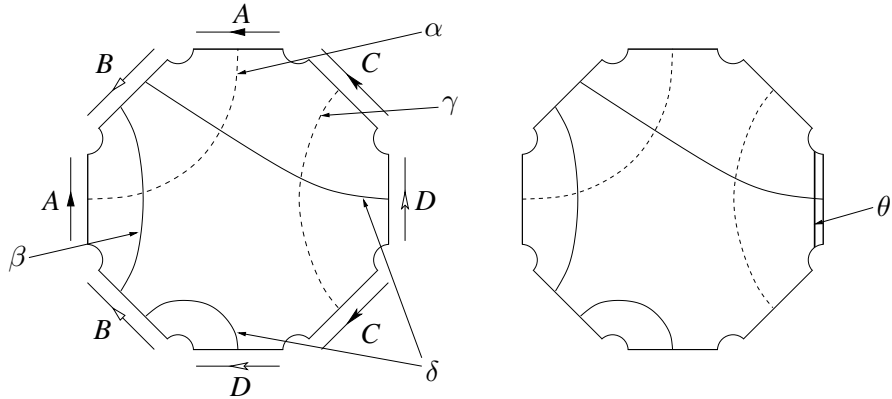

In the picture there are represented four further curves:
$\alpha$, $\beta$, $\gamma$ and $\delta$. Defining
\begin{align}
\cal C=&\{\,\beta,\delta\,\},\notag\\
\cal D=&\{\,\alpha,\gamma\,\},\notag
\end{align}
it is easy to check that $(\cal C,\cal D)$ is a Penner pair in
$S$. The automorphism $\varphi\co S\to S$ defined by
$$
\varphi=T_\beta^-\circ T_\delta^-\circ T_\alpha^+\circ T_\gamma^+
$$
is, therefore, a Penner automorphism subordinate to the pair
$(\cal C,\cal D)$.

The pair $(\cal C,\cal D)$ admits dual arcs. The picture shows
one, labelled as $\theta$. We consider the corresponding disc
$\Delta_{\theta}$. \fullref{F:twists in product} shows
$S_0=S\times\{0\}$, $S_1=S\times\{1\}\subseteq\partial H$ and how
$\partial\Delta_{\theta}$ intersects them\footnote{In fact, one
can picture the whole $\partial\Delta_\theta\cap\partial H$ in the
figure. The only portion of $\partial H$ not represented is the
vertical annulus $\partial S\times I$, which
$\partial\Delta_\theta$ intersects in four vertical arcs.}.

\begin{figure}[ht!]
\labellist\small
\hair1pt
\pinlabel {$S\times\{0\}$} at 100 358
\pinlabel {$S\times\{1\}$} at 570 358
\pinlabel {$\partial\!\Delta_\theta$} <0.5pt,-1.5pt> at 417 66
\endlabellist
\centerline{\includegraphics[scale=0.35]{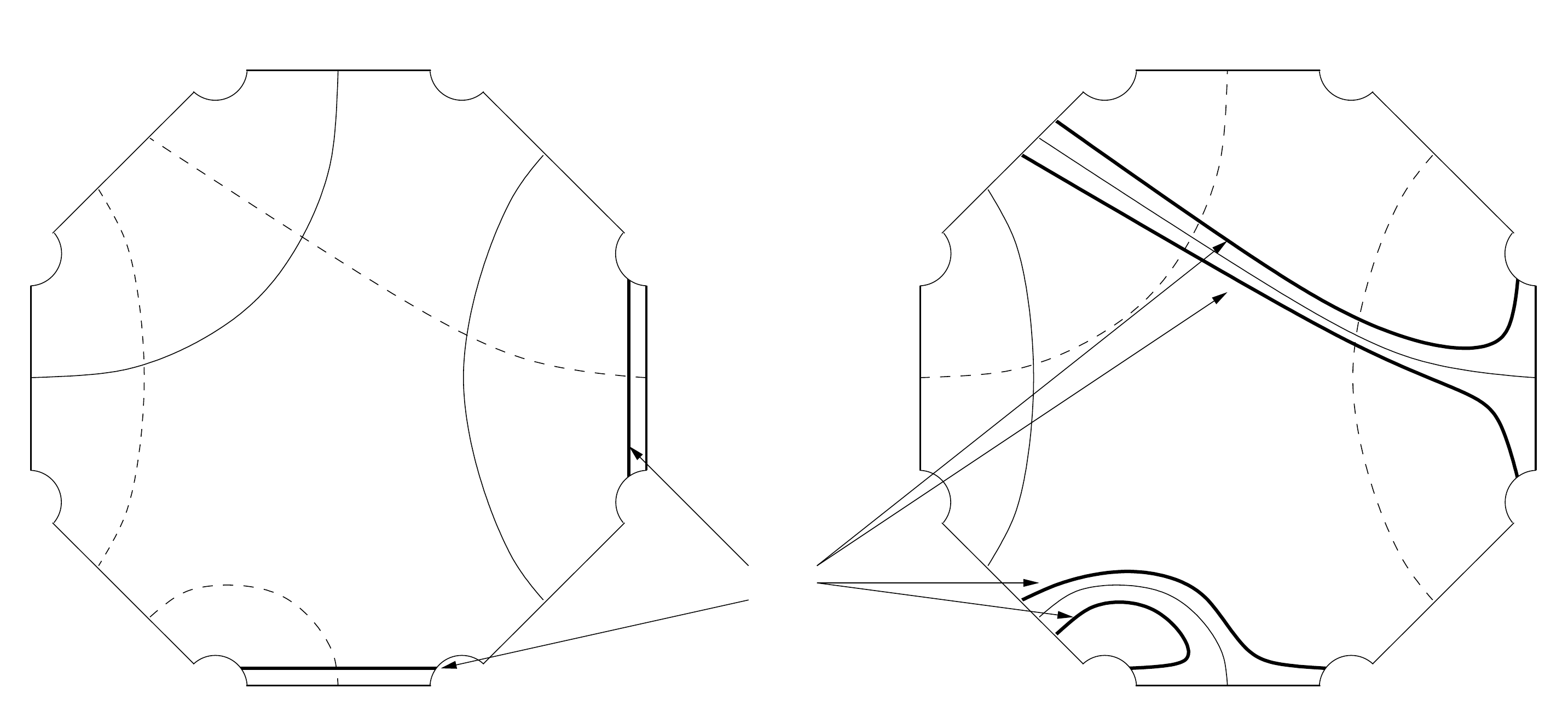}}
\caption{The curve $\partial\Delta_{\theta}$ in $\partial H$}
\label{F:twists in product}
\end{figure}

\fullref{F:twists in product} also shows the pair
$(\mathcal{Q},\mathcal{R})$ obtained by \fullref{L:twisting
disc}: $\mathcal{Q}$ consists of the ``solid'' lines, including
$\partial\Delta_{\theta}$, while the ``dashed'' lines form
$\mathcal{R}$.

Let $\widehat\varphi\co H\to H$ be the lift of $\varphi$ to $H$. By
\fullref{T:method2}
$$
\widehat\varphi\circ T^-_{\Delta_{\theta}} \co H\to H
$$
is an irreducible automorphism, where $T^-_{\Delta_{\theta}}$ is
the left twist along $\Delta_{\theta}$.
\end{example}

%%%%%%%%%%%%%%%%%%%%%%%%%%%%%%%%%%%%%%%%%%%%%%%%%%%%%%%%%%%%%%%%%%%%%%%%%%%%%%%%%%%%%%%%%%
%%%%%%%%%%%%%%%%%%%%%%%%%%%%%%%%%%%%%%%%%%%%%%%%%%%%%%%%%%%%%%%%%%%%%%%%%%%%%%%%%%%%%%%%%%

\subsection{The incompressibility condition} \label{SS:develop
second}

\begin{example}\label{E:incompressible}
We continue working with \fullref{E:method}. Now we will
determine certain $f$--invariant laminations and estimate their
corresponding growth rate $\lambda$. Recall \fullref{SS:background},
especially \fullref{T:generic}, \fullref{Rs:laminations} and
\fullref{T:incompressibility}. We
refer the reader to \cite{UO:Autos} for details on the
constructions.

Recall the oriented surface $S$, a once punctured genus two
surface. Here it will be convenient to regard $S\times I$ as
$H_0\subseteq H$. We now choose a complete system of discs
$\cal{E}_0$ in $H_0$, as follows. Consider the labelled arcs $A$,
$B$, $C$ and $D\subseteq S$ as in the figure (the sides of the
octagon), and construct the discs $A\times I$, $B\times I$,
$C\times I$, $D\times I\subseteq H_0$. Abusing notation, we use
the same labels $A$, $B$, $C$ and $D$ respectively to represent
these discs. Let $\cal{E}_0=\{A, B, C, D\}$. From $\cal{E}_0$ we
can consider the dual graph $\Gamma$. In fact, we will regard
$\Gamma$ as an oriented labeled graph (represented as a spine in
\fullref{F:image in product}, on the left). To avoid
ambiguities we require that $\Gamma\subseteq
S\times\{\frac{1}{2}\}$. We now regard $H_0$ as a neighborhood of
$\Gamma$ and consider the corresponding handle decomposition
$\cal{H}_0$ in $H_0$.

Together with the handle decomposition we will now choose a
representative in the class of $f$ (which we also label as $f$)
and study the associated laminations. These are determined by how
the handles of $H_0$ intersect the handles of $H_1=f(H_0)$. But it
is equivalent to consider $H_{-1}=f^{-1}(H_0)\subseteq H_0$ (apply
the diffeomorphism $f^{-1}$), which is easier to picture.

Regarding $H_0$ as a neighborhood of $\Gamma$, we consider
$f^{-1}(\Gamma)~=T_{\Delta_\theta}^+\circ(\widehat\varphi)^{-1}(\Gamma)\subseteq
H_0$ (the disc $\Delta_\theta$ and the automorphism $\widehat\varphi$
are defined in \fullref{E:method}). \fullref{F:image in
product}, right, shows this image $f^{-1}(\Gamma)$, determining
how $\Gamma$, and hence $H_0$, should be pictured in $H_1$.

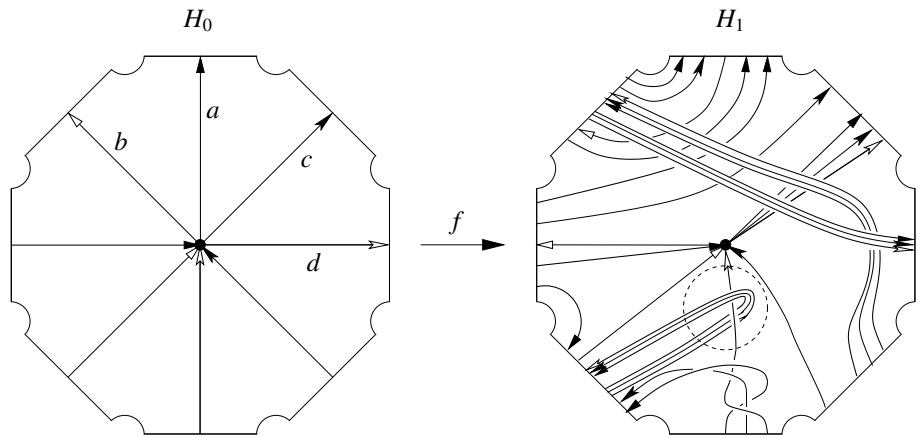
\begin{figure}[ht!]
\labellist\small
\pinlabel {$a$} [l] at 177 292
\pinlabel {$b$} [bl] at 100 255
\pinlabel {$c$} [tl] at 260 260
\pinlabel {$d$} [t] at 276 178
\pinlabel {$f$} [b] at 400 180
\pinlabel {$H_0$} [b] at 177 355
\pinlabel {$H_1$} [b] at 634 355
\endlabellist
\centerline{\includegraphics[scale=0.45]{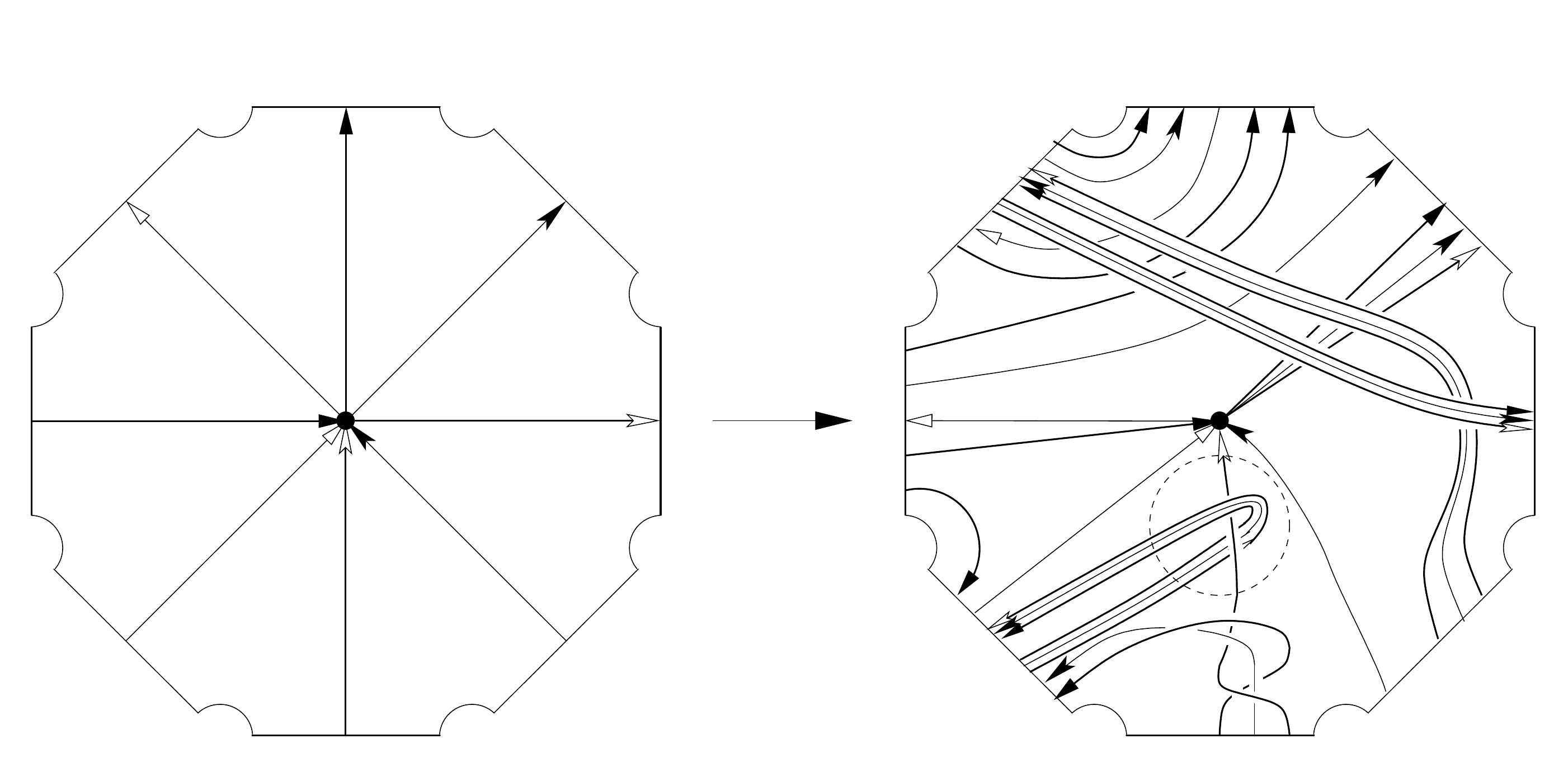}}
\caption{The graph $\Gamma\subseteq H_0$ and $\Gamma\subseteq H_1$}
\label{F:image in product}
\end{figure}

This intersection pattern determines an $f$--invariant measured
lamination $\Lambda$ with full support provided that the incidence
matrix corresponding to $\cal{E}_0$ is irreducible. We verify this
using the transpose of that matrix, which is the incidence matrix
$M(\cal{H}_0)$ for the $1$--handles of $\cal{H}_0$. To do that
consider the handle decomposition $\cal{H}_1$ of $H_1$ induced by
$\cal{H}_0$ through $f$. The incidence matrix
$M=M(\cal{H}_0)=\{\,m_{ij}\,\}$ is given by $m_{ij}=|f(e_i)\cap
e_j|$, where $f(e_i)$ is a $1$--handle of $\cal{H}_1$ and $e_j$ is
a $1$--handle of $\cal{H}_0$. Therefore
$$
M=M(\cal{H}_0)=\left(
\begin{matrix}
3&1&1&0\\
4&1&3&2\\
1&0&2&1\\
1&0&1&1
\end{matrix}
    \right),\notag
$$
which is indeed irreducible (all entries of $M^2$ are strictly
positive) with a Perron--Frobenius eigenvalue, to three decimal places, of
$$
\lambda=\lambda(\cal{H}_0)\approx 4.987.
$$
We now ask the question of whether $\lambda$ is the minimal growth
in the isotopy class of $f$ or not. The ``no back-tracking''
sufficient condition does not apply to this case (for instance,
there is an unremovable ``back-tracking'', linked on the circled
part of \fullref{F:image in product}, right). Oertel's
incompressibility property, a necessary condition, holds here, as
methodic --- though tedious --- computation reveals\footnote{One
can carry on this computation in a manner not unlike those in the papers 
\cite{BH:Surfaces,BH:Tracks} by Bestvina and Handel. Here one just has to be
careful with certain back-trackings, which are allowed because of
linkings (eg the circled section in \fullref{F:image in
product}).}.

We see then that, up to this point, we have no criterion to decide
whether $\lambda\approx 4.987$ is minimal or not. In the next
section, \fullref{S:tightness}, a simple argument shows that
this growth rate is not minimal. This will be related to the fact
that this lamination is not ``tight''.

\end{example}

%%%%%%%%%%%%%%%%%%%%%%%%%%%%%%%%%%%%%%%%%%%%%%%%%%%%%%%%%%%%%%%%%%%%%%%%%%%%%%%%%%%%%%%%%%
%%%%%%%%%%%%%%%%%%%%%%%%%%%%%%%%%%%%%%%%%%%%%%%%%%%%%%%%%%%%%%%%%%%%%%%%%%%%%%%%%%%%%%%%%%
%%%%%%%%%%%%%%%%%%%%%%%%%%%%%%%%%%%%%%%%%%%%%%%%%%%%%%%%%%%%%%%%%%%%%%%%%%%%%%%%%%%%%%%%%%
%%%%%%%%%%%%%%%%%%%%%%%%%%%%%%%%%%%%%%%%%%%%%%%%%%%%%%%%%%%%%%%%%%%%%%%%%%%%%%%%%%%%%%%%%%

\section{Tightness}\label{S:tightness}

From now on we assume that $f\co H\to H$ is an irreducible
automorphism of a handlebody $H$. We consider a handle
decomposition $\cal{H}_0$ of $H_0$ and the corresponding: disc
system $\cal{E}_0$, $f$--invariant measured laminations
$(\Lambda,\mu)$, $(\Omega,\nu)$ and growth rate $\lambda$. The
one--dimensional lamination $\Omega$ and its measure $\nu$ will
play important roles throughout this section. Given an immersed
surface $F\subseteq H$ transverse to $\Omega$ we will denote
$\int_F\nu$ by either $\nu(F)$ or $\inter{F}{\Omega}{\nu}$. The
advantage of the first notation is in its simplicity and will be
preferred whenever there are no ambiguities. The advantage of the
second is that it emphasizes the object $\Omega$ supporting $\nu$,
which will be convenient in certain contexts. We call
$\nu(F)=\inter{F}{\Omega}{\nu}$ the \emph{weighted intersection}
(or just the \emph{intersection}) of $F$ with $(\Omega,\nu)$.

Recall that the goal is to characterize minimal growth.

\subsection{Tightening discs}\label{SS:tightening discs}

%%%%%%%%%%%%%%%%%%%%%%%%%%%%%%%%%%%%%%%%%%%%%%%%
%
%   Motivating example
%
%%%%%%%%%%%%%%%%%%%%%%%%%%%%%%%%%%%%%%%%%%%%%%%%

In \fullref{E:incompressible} we left unproven the claim that
Oertel's incompressibility property does not imply minimal
$\lambda$. Recall from \fullref{T:incompressibility} the
incompressibility property, which can also be stated as follows.
The two--dimensional lamination $\Lambda$ has the property if for
any $n>0$ the leaves of $\Lambda\cap(H_n-\ring{H}_0)$ (which are
properly embedded planar surfaces with a boundary component in
$\partial H_n$ and the others in $\partial H_0$) are
incompressible in $H_n-\ring{H}_0$ (see the paper \cite{UO:Autos} by
Oertel).

Let then $f\co H\to H$ be an irreducible automorphism, assume that
an invariant $\Lambda$ has the incompressibility property and
consider the associated growth rate $\lambda$. The next schematic
example suggests a reason for the fact that this incompressibility
does not imply minimality of $\lambda$. A handle decomposition
$\cal{H}_0$ of $H_0$ determines, through $f^i$, a handle
decomposition $\cal{H}_i$ of $H_i=f^i(H_0)$, $i\in\mathbb{Z}$. The
corresponding incidence matrix $M(\cal{H})$ is assumed to be
irreducible, a condition required for the construction of the
laminations (see the paragraph preceding \fullref{P:subinvariance}).

\begin{example}\label{E:loose}
Let $V$ be a $0$--handle of $\cal{H}_1$ and suppose that $H_0\cap
V$ is as in \fullref{F:loose picture}(a), with $V$ intersecting
$1$--handles $e_p$ and $e_q$ of $\cal{H}_0$, and $f(E)$ the image
of a disc $E\in\cal E_0=\{\,E_1,\dots, E_k\,\}$.

\begin{figure}[ht!]
\labellist\small
%\psfrag{ei}{$e_p$}
%\psfrag{ej}{$e_q$}
%\psfrag{V}{$V$}
%\psfrag{f(E)}{$f(E)$}
\pinlabel {$e_p$} [br] at 140 174
\pinlabel {$e_p$} [tl] at 612 153
\pinlabel {$e_q$} [tl] at 230 302
\pinlabel {$e_q$} [tl] at 720 302
\pinlabel {$V$} at 130 300
\pinlabel {$V$} at 615 300
\pinlabel {$f(E)$} at 24 145
\pinlabel {$f(E)$} [l] at 602 225
\pinlabel {(a)} [b] at 285 0
\pinlabel {(b)} [b] at 770 0
\endlabellist
\centerline{\includegraphics[scale=0.35]{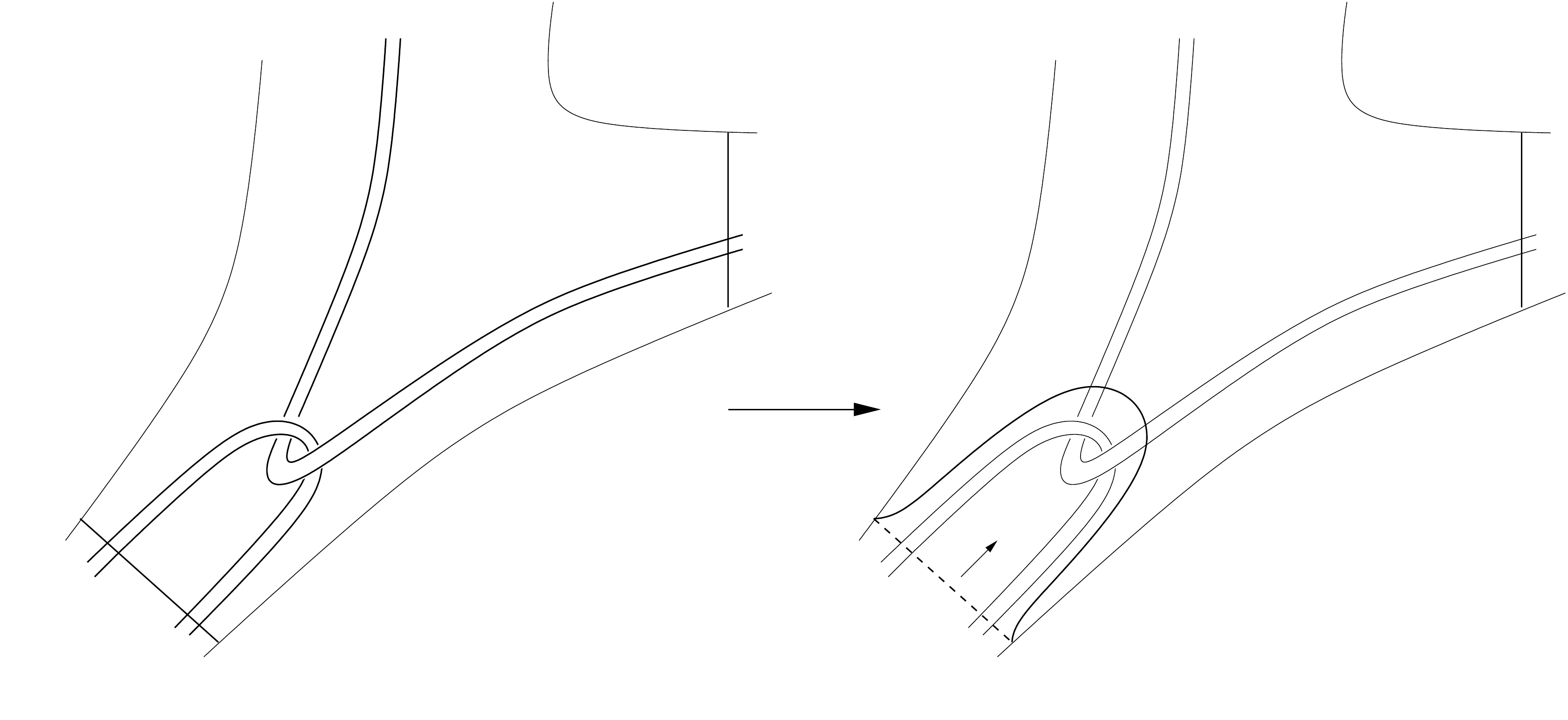}}
\caption{(a)\qua A disc $f(E)$ dual to a $1$--handle of $H_1$,
intersecting the $1$--handle $e_i$ of $\cal{H}_0$;\qua (b)\qua the move
that reduces $\lambda$}\label{F:loose picture}
\end{figure}
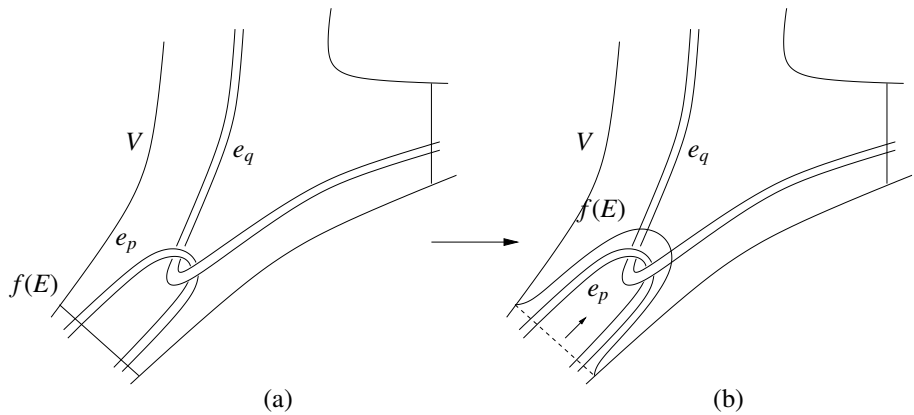

Now consider $M=M(\cal{H}_0)=\{\,m_{ij}\,\}$, the incidence matrix
for the $1$--handles of $\cal{H}_0$: $m_{ij}=|f({e}_i)\cap {e}_j|$
counts how many times the $1$--handle $e_j$ (of $\cal{H}_0$)
crosses the $1$--handle $f(e_i)$ of $\cal{H}_1$. Recall that we
are assuming that $M$ is irreducible. Its Perron--Frobenius
eigenvalue is then $\lambda$, the growth rate of $f$ (with respect
to $\cal{H}_0$).

Now note that $E_i$ is transverse to $\Omega$, so it makes sense
to consider
$$
\widehat\nu_i=\nu(E_i),
$$
where we recall that $\nu$ is the transverse measure on $\Omega$.
That determines a vector $\widehat\nu=(\widehat\nu_1,\dots,\widehat\nu_k)$,
which is precisely a Perron--Frobenius eigenvector of $M$:
$$
M\widehat\nu=\lambda\widehat\nu.
$$
Suppose that $\widehat\nu_q<\widehat\nu_p$. We can isotope $f$ to replace
intersections of $f(E)$ with $e_p$ by intersections with $e_q$
(see \fullref{F:loose picture}(b)).

This operation does not change the handle decomposition. Now the
new incidence matrix $M'=\{\,m'_{ij}\,\}$ is given by
\begin{align*}
m'_{ij}& =m_{ij} & \text{if $ij\neq 1p$, $1q$;}\\
m'_{1p}& =m_{1p}-2;&\\
m'_{1q}& =m_{1q}+2. &
\end{align*}
Suppose that $M'$ is irreducible. In this case recall that $\widehat
\nu_q<\widehat \nu_p$ and consider $M'\widehat \nu$:
\begin{align*}
(M'\widehat\nu)_i & =\lambda\widehat\nu_i & \text{if $i\neq 1$;} \\
(M'\widehat\nu)_1 &
=\lambda\widehat\nu_1-2\widehat\nu_p+2\widehat\nu_q<\lambda\widehat\nu_1. &
\end{align*}
By \fullref{P:subinvariance}
$$\lambda(M')<\lambda=\lambda(M),$$
therefore the isotopy reduces the growth rate.
\end{example}

%%%%%%%%%%%%%%%%%%%%%%%%%%%%%%%%%%%%%%
%
%   First definitions
%
%%%%%%%%%%%%%%%%%%%%%%%%%%%%%%%%%%%%%%%%%%%%%%%%%%%%%%%%%%%%%%%

A situation like the one described in the example above indeed
happens, see \fullref{E:loose2}. It not only shows that
Oertel's incompressibility property does not imply minimality of
the growth rate but also suggests that the weighted intersection
$\inter {E}{\Omega}{\nu}=\nu(E)$ (where $E\in\cal E_0$) should be
relevant in the search for the minimal growth. We introduce, then,
the following definition:

\begin{definition}\label{D:tight}
Let $(\Delta,\partial\Delta)\to(H,\Lambda)$ be an embedded disc
transverse to $\Omega$. Consider $\Delta'\subseteq\Lambda$ such
that $\partial\Delta'=\partial\Delta$.  We say that $\Delta$ is a
\emph{\tightening\ disc} for the triple $(\Lambda,\Omega,\nu)$ if
$$\nu(\Delta)<\nu(\Delta').$$
The triple $(\Lambda,\Omega,\nu)$ is said \emph{\tight} if there is
no \tightening\ disc. We will often abuse notation and say that
$\Lambda$ is tight or not, leaving $(\Omega,\nu)$ implicit.
Accordingly, we may say that a tightening disc for
$(\Lambda,\Omega,\nu)$ is a tightening disc for $\Lambda$ only.
\end{definition}
\begin{obs}\label{R:tightening disc}

The requirement that the \tightening\ disc $\Delta$ is transverse
to $\Omega$ implies that it does not intersect the singular set
$S(\Omega)$.

We also note that $\Lambda$ being \tight\ implies that it has the
incompressibility property: a compressing disc for $\Lambda-\ring
H_0$ is a \tightening\ disc.
\end{obs}

Now we can say that the original lamination in
\fullref{E:loose} is not tight, with a \tightening\ disc
represented in \fullref{F:loose picture}(b). As previously
mentioned, that is a hypothetical situation. The following is a
specific example.

\begin{example}\label{E:loose2}
We refer to \fullref{E:incompressible} and consider the
automorphism $f\co H\to H$ and the handle decomposition
$\cal{H}_0$ of $H_0$ defined then. We consider the disc
$f(B)\subseteq H_1$, the co-core of a handle of $\cal{H}_1$.

\begin{figure}[ht!]
\labellist\small
%\psfrag{a}{$a$}
%\psfrag{b}{$b$}
%\psfrag{c}{$c$}
%\psfrag{d}{$d$}
%\psfrag{alpha}{$A$}
%\psfrag{beta}{$B$}
%\psfrag{gamma}{$C$}
%\psfrag{delta}{$D$}
%\psfrag{f(beta)}{$f(B)$}
%\psfrag{f}{$f$}
%\psfrag{H0}{$H_0$}
%\psfrag{H1}{$H_1$}
%\psfrag{D}{$\Delta$}
%\psfrag{D'}{$\Delta'$}
\pinlabel {$a$} [l] at 186 300
\pinlabel {$b$} [bl] at 109 263
\pinlabel {$c$} [tl] at 268 268
\pinlabel {$d$} [t] at 280 186
\pinlabel {$A$} [b] at 186 370
\pinlabel {$B$} [br] at 60 318
\pinlabel {$B$} [tr] at 60 60
\pinlabel {$C$} [bl] at 310 310
\pinlabel {$D$} [l] at 364 160
\pinlabel {$f$} [b] at 400 185
\pinlabel {$f(B)$} [tr] at 500 95
\pinlabel {$H_0$} [b] at 55 370
\pinlabel {$H_1$} [b] at 740 370
\pinlabel {$\Delta$} [l] at 792 56
\pinlabel {$\Delta'$} [tr] at 477 40
\endlabellist
\centerline{\includegraphics[scale=0.45]{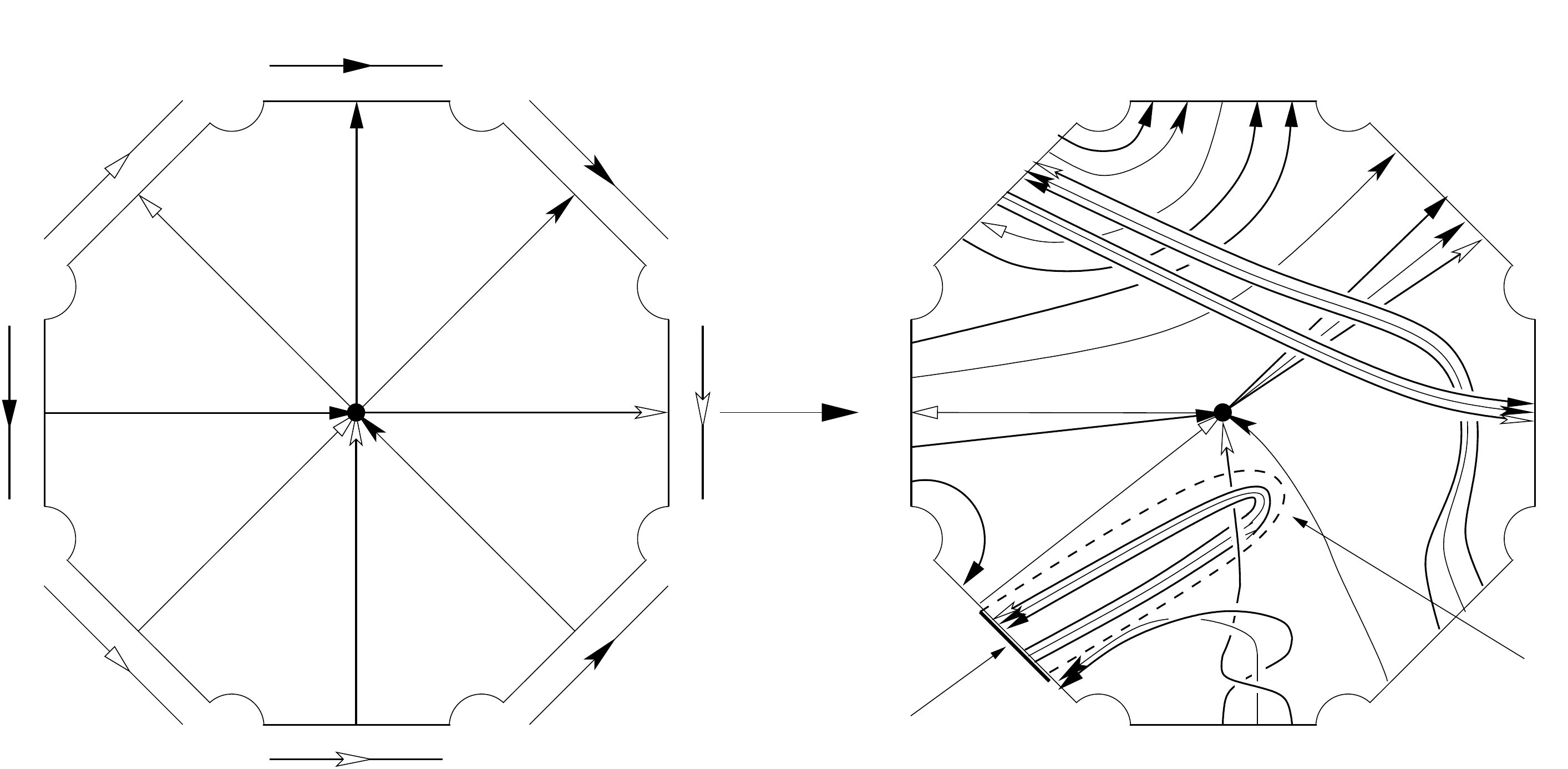}}
\caption{A \tightening\ disc $\Delta$ in $H_1$}\label{F:tightening disc}
\end{figure}

One can see in \fullref{F:tightening disc} a \tightening\ disc
$\Delta$ (represented at the right by a dashed line), with its
boundary in a leaf of $\Lambda\cap H_1$ parallel to $f(B)$.
Indeed, let $\Delta'\subseteq f(B)$ be the disc such that
$\partial\Delta'=\partial\Delta$ (represented by a thick line). It
is easy to check that
$$
\nu(\Delta)=2\nu(D)< 2(\nu(A)+\nu(C)+\nu(D))=\nu(\Delta'),
$$
for $\nu$ has full support on $\Omega$. We now can change $f$
through an isotopy taking $\Delta'$ to $\Delta$.
\fullref{F:pulling tight} shows the result of such an isotopy.

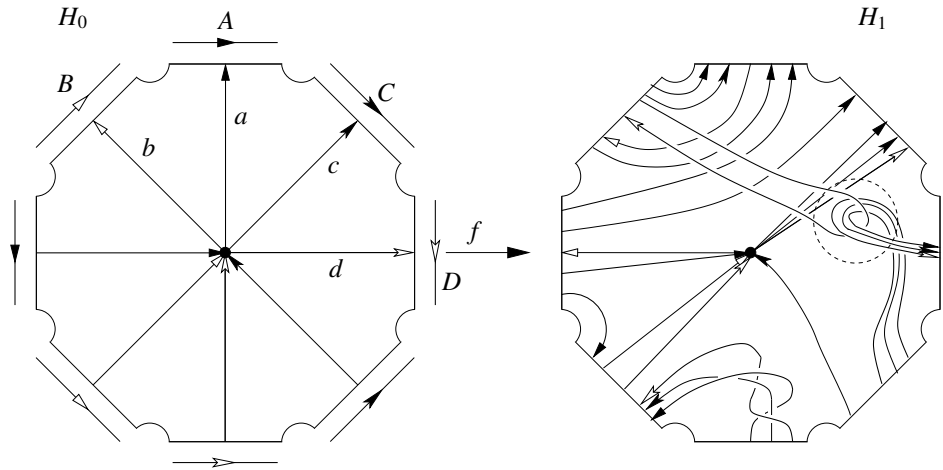
\begin{figure}[ht!]
\labellist\small
%\psfrag{a}{$a$}
%\psfrag{b}{$b$}
%\psfrag{c}{$c$}
%\psfrag{d}{$d$}
%\psfrag{alpha}{$A$}
%\psfrag{beta}{$B$}
%\psfrag{gamma}{$C$}
%\psfrag{delta}{$D$}
%\psfrag{f}{$f$}
%\psfrag{H0}{$H_0$}
%\psfrag{H1}{$H_1$}
\pinlabel {$a$} [l] at 186 300
\pinlabel {$b$} [bl] at 109 263
\pinlabel {$c$} [tl] at 268 268
\pinlabel {$d$} [t] at 280 186
\pinlabel {$A$} [b] at 186 370
\pinlabel {$B$} [br] at 60 318
\pinlabel {$C$} [bl] at 310 310
\pinlabel {$D$} [l] at 364 160
\pinlabel {$f$} [b] at 400 185
\pinlabel {$H_0$} [b] at 55 370
\pinlabel {$H_1$} [b] at 740 370
\endlabellist
\centerline{\includegraphics[scale=0.45]{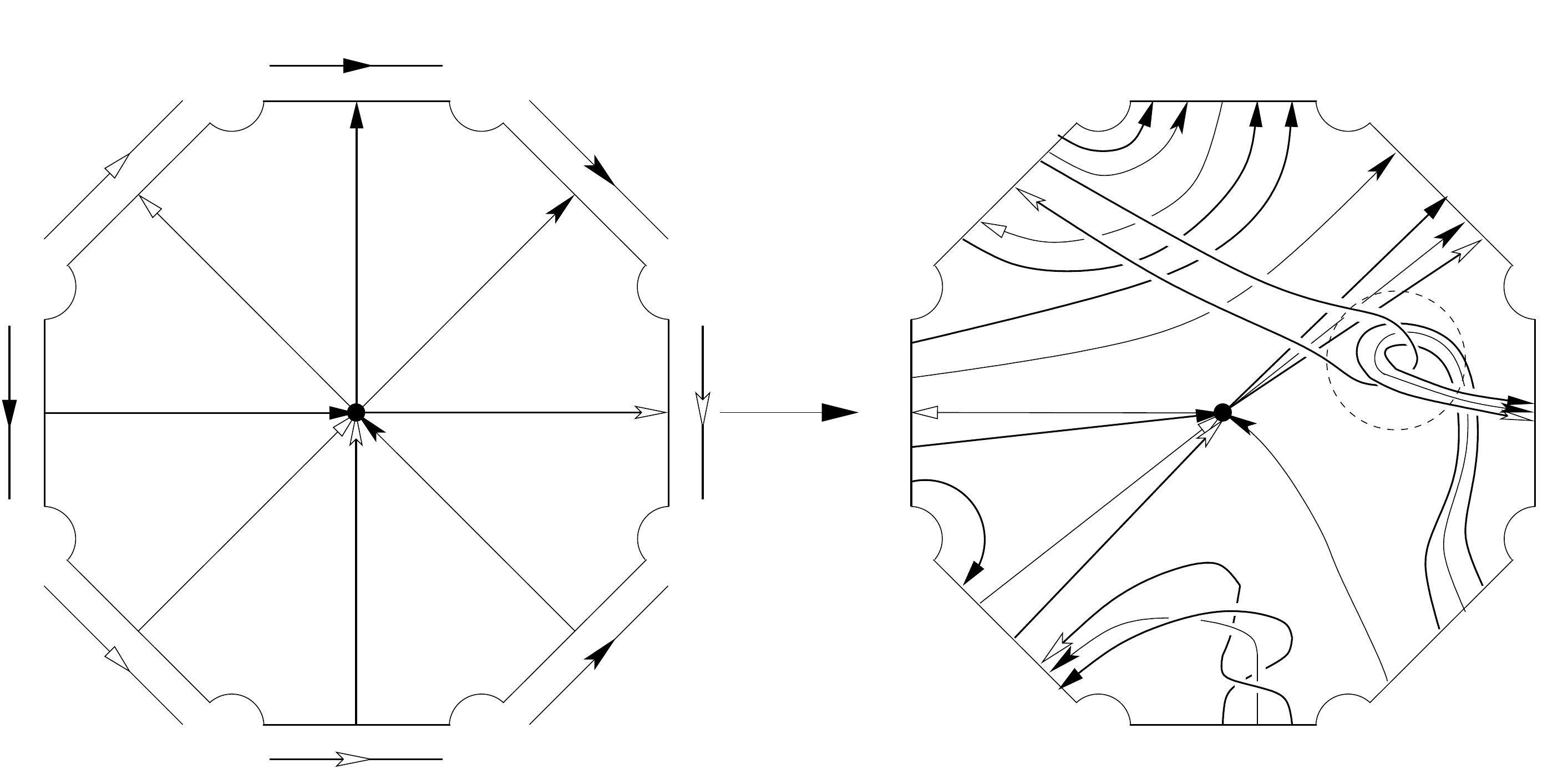}}
\caption{After the growth-reducing isotopy}
\label{F:pulling tight}
\end{figure}

We verify that the new incidence matrix is:
$$
M(\cal{H})=\left(
\begin{matrix}
3&1&1&0\\
2&1&1&2\\
1&0&2&1\\
1&0&1&1
\end{matrix}
    \right)\notag
$$
which is irreducible with Perron--Frobenius eigenvalue
$\lambda(\cal{H})\approx 4.542$ (precise up to three decimals),
showing that the previous lamination did not have minimal growth
rate.
\end{example}

We are interested in the problem of characterizing minimal growth
rate. Considering these examples one should expect \tightness\ to
play a role in the solution.

\begin{conjecture}\label{Cj:equivalence}
The growth rate is minimal if and only if the lamination $\Lambda$
is \tight.
\end{conjecture}

At this point we note that tightness is strictly stronger than
Oertel's incompressibility --- which is too weak --- and strictly
weaker than ``no back-tracking'' --- which is too strong.

In the direction of proving the conjecture we will prove the following
theorem.

\begin{theorem}\label{T:minimality}
If $\Lambda$ is tight then $\lambda$ is minimal.
\end{theorem}

Concerning the converse we will prove the following proposition.

\begin{proposition}\label{P:valence}
Suppose that all $0$--handles of the handle decomposition
$\cal{H}$ have \emph{valence} $2$ or $3$. If $\lambda$ is minimal
then $\Lambda$ is tight.
\end{proposition}

By \emph{valence} of a $0$--handle of $\cal{H}_0$ we mean the
number of ends of $1$--handles that are attached to it. In other
words, consider the graph $\Gamma$ corresponding to $\cal{H}_0$.
The valence of the vertex corresponding to the $0$--handle is its
valence. As an application of the technical proposition above we
will get the following corollary.

\begin{corollary}\label{C:true genus 2}
\fullref{Cj:equivalence} is true if the handlebody has
genus 2.
\end{corollary}

The proofs will be given in the following subsections.

%%%%%%%%%%%%%%%%%%%%%%%%%%%%%%%%%%%%%%%%%%%%%%%%%%%%%%%%%%%%%%%%%%%%%%%%%%%%%%%%%%%%%%%%%%
%%%%%%%%%%%%%%%%%%%%%%%%%%%%%%%%%%%%%%%%%%%%%%%%%%%%%%%%%%%%%%%%%%%%%%%%%%%%%%%%%%%%%%%%%%

\subsection{Strong tightening discs}\label{SS:strong}

In this subsection we introduce some technical constructions and
results.

The definition of \tightening\ disc that we gave was quite
general. We will work with tightening discs having some special
properties. These will be called \emph{strong tightening disc} (see
\fullref{D:strong} below). We will show that there is no
loss of generality in working with them (\fullref{P:equivalence}).

Consider the handle decomposition $\cal{H}_0$ of $H_0$, with
1--handles $e_1$, $\dots$, $e_k$. We will give weights to these
$1$--handles, ie assign a positive number $v_j$ to each
$1$--handle $e_j$. More precisely, given a positive vector
$v=(v_1,\dots,v_k)$, ($v_j>0$), we say that the pair
$(\cal{H}_0,v)$ is a \emph{weighted handle decomposition of $H_0$}.
In this case we say that $v$ is a \emph{system of weights in
$\cal{H}_0$}. We extend these notions for any level $i$ through
$f^i$, so it makes sense to refer to a weighted handle
decomposition $(\cal{H}_i,v)$ for any $H_i=f^i(H_0)$.

Consider a weighted handle decomposition $(\cal{H}_i,v)$ for $H_i$
and let $S$ be an embedded surface intersecting $H_i$ just in its
1--handles. Suppose further that $S\cap H_i$ consists of discs
dual to the 1--handles $e_j$ of $\cal{H}_i$. We define
$$\inter S {\cal{H}_i} v=\sum_{1\leq j\leq k} |S\cap e_j|\cdot v_j.$$
If $E_j$ is a disc dual to the $1$--handle $e_j$ of
$\mathcal{H}_0$ we consider
$$\widehat\nu_j=\inter{E_j}{\Omega}{\nu}=\nu(E_j).$$
Note that $\widehat\nu_j$ does not depend on the choice of dual disc
$E_j$. Also, the vector $\widehat\nu=(\widehat\nu_1,\dots,\widehat\nu_k)$ is a
Perron--Frobenius eigenvector of the incidence matrix associated
to $\cal{H}_0$ and $f$.

\begin{definition}\label{D:standard}
Consider $\widehat\nu$ built above. We regard it as a system of
weights in $\cal{H}_0$. We call $(\cal{H}_0,\widehat\nu)$ the \emph{standard
weighted handle decomposition of $H_0$}.
\end{definition}

\begin{obs}\label{R:standard}
It is clear by the construction above that
\begin{equation}\notag
\inter{E_j}{\cal{H}_0}{\widehat\nu}=\inter{E_j}{\Omega}{\nu}.
\end{equation}
Hence, for a general surface $S$ intersecting $H_0$ in dual discs,
\begin{equation}\label{Eq:inter}
\inter{S}{\cal{H}_0}{\widehat\nu}=\inter{S}{\Omega}{\nu}.
\end{equation}
\end{obs}

In the definition below $(\cal{H}_0,\widehat\nu)$ is standard.

\begin{definition}\label{D:strong}
Let $(\Delta,\partial\Delta)\subseteq (H,\Lambda)$ be an embedded
disc and $\Delta'\subseteq\Lambda$ be such that
$\partial\Delta'=\partial\Delta$. We say that $\Delta$ is a \emph{strong
tightening disc} if there exists $n$ such that
\begin{enumerate}
  \item $\Delta\subseteq H_n$ and $\partial\Delta\subseteq\partial
  H_n$;\label{D:strong-p1}
  \item for any $0\leq i\leq n$, $\Delta\cap H_i$ consists of essential
  discs (in $H_i$) and $\Delta\cap H_0$ consists, moreover, of discs dual
  to the $1$--handles;\label{D:strong-p2}
  \item $\Delta\cap\Lambda=\partial\Delta$;\label{D:strong-p3}
  \item $\inter{\Delta}{\cal{H}_0}{\widehat\nu}<
         \inter{\Delta'}{\cal{H}_0}{\widehat\nu}$.\label{D:strong-p4}
\end{enumerate}
\end{definition}

\begin{proposition}\label{P:equivalence}
There exists a \tightening\ disc if and only if there exists a
strong tightening disc.
\end{proposition}
\begin{proof}
One direction is immediate: a strong tightening disc is a
\tightening\ disc by property 4 and the equation \eqref{Eq:inter}
above.

We prove the other direction. Let $\Delta$ be a tightening disc.
We will build another one with the required properties. These
properties will be realized progressively (not necessarily in the
order specified in \fullref{D:strong} above). For
simplicity of notation we will also label these intermediate discs
as $\Delta$.

\textbf{Part 1: property \fullref{D:strong-p3}}\qua
Let $n$ be such that $\Delta\subseteq H_n$. We first note that
$\partial\Delta$ is contained in a leaf of $\Lambda\cap H_n$,
which consists of discs dual to the 1--handles of $H_n$ (see
\fullref{Rs:laminations}). In particular $\partial\Delta$ is
contained in a 1--handle $e_l$ and fix
$E_l\supseteq\partial\Delta$ the corresponding dual disc. In each
other 1--handle $e_i$ of $H_n$ choose an arbitrary dual disc
$E_i\subseteq H_n$. We assume that $\Delta$ is transverse to
$\bigcup_i E_i$, including at $\partial\Delta\subseteq E_j$.

The main goal is to
reduce the complexity $\bigl|\Delta\cap \bigcup_i E_i\bigr|$ by
performing surgeries and isotopies which preserve the
property of being a tightening disc, to eventually yield
$\Delta\cap\bigcup_i E_i=\partial\Delta$.

Consider $\Delta\cap\bigcup_i E_i$ (which consists just of closed
curves) and choose a curve $\gamma\subseteq E_j$ that is innermost
in some $E_j$. Let $D\subseteq E_j$ and $D'\subseteq\Delta$ be the
discs bounded by $\gamma$. There are two cases to consider:

\textbf{Case 1: $\nu(D)\leq\nu(D')$}\qua
    In this case we perform a surgery in $\Delta$, replacing
    $D'\subseteq\Delta$ by $D$ and pushing it a bit away from $E_j$
    (this pushing
    should be vertical, that is, along the $I$--fibers of the product
    structure $D^2\times I$ of the $1$--handle $e_j$, and ``to the
    side opposite to'' $D'$).  Since $\Omega\cap e_j$ consists of
    $I$--fibers, which are preserved by the ``pushing'' move, the
    process does not increase
    $\nu(\Delta)$. It clearly reduces $\bigl|\Delta\cap\bigcup_i E_i\bigr|$.

\textbf{Case 2: $\nu(D)>\nu(D')$}\qua
    Here $D'$ is a tightening disc by definition. If
    $\gamma\neq\partial\Delta$ then $\bigl|D'\cap\bigcup_i
    E_i\bigr|<\bigl|\Delta\cap\bigcup_i
    E_i\bigr|$ and we replace $\Delta$ with $D'$, reducing complexity.
    If $\gamma=\partial\Delta\subseteq E_l$ then there are several
    possibilities:
    \begin{itemize}

        \item There exists another curve
        $\gamma'\subseteq\bigl(\Delta\cap\bigcup_i E_i\bigr)$,
        $\gamma'\neq\gamma$,
        innermost in some $E_j$. In this case $\gamma'\neq\partial\Delta$
        and we apply the procedure described above for $\gamma'$,
        reducing $\bigl|\Delta\cap\bigcup_i E_i\bigr|$.

        \item $\Delta\cap\bigcup_i E_i=\gamma=\partial\Delta$ and we have
        achieved what was desired.

    \item $\bigl|\Delta\cap\bigcup_i E_i\bigr|\geq 2$ and no curve
    $\gamma'\subseteq\bigl(\Delta\cap\bigcup_i E_i\bigr)$,
    $\gamma'\neq\gamma=\partial\Delta$
    is innermost in the dual disc $E_j$ that contains it.
    Therefore $\partial\Delta\subseteq E_j=E_l$ and, if
    $i\neq r$, $\Delta\cap E_i=\emptyset$. Let then $\gamma'\subseteq
    \Delta\cap E_l$ be \emph{second innermost} in $E_l$, in the sense
    that the (interior of the) disc $D''\subseteq E_l$ it
    bounds contains just innermost curves. Since we are assuming
    that $\partial\Delta$ is the single innermost curve, $D''\cap
    \Delta=\gamma'\cup\partial\Delta$. On the other hand,
    $\gamma'\subseteq
    \Delta$ bounds a disc $\Delta''\subseteq\Delta$. We shall prove
    that $\Delta''$
    is a \tightening\ disc. Indeed recall that $\Delta''\subseteq
    \Delta$ and $\Delta'\subseteq D''$ and, therefore
    $$\nu(\Delta'')\leq\nu(\Delta)\quad {\rm and}\quad
      \nu(\Delta')\leq\nu(D'').$$
    But $\Delta$ is a \tightening\ disc, ie $\nu(\Delta)<\nu(\Delta')$.
    Combining these inequalities one gets that
    $$\nu(\Delta'')<\nu(D''),$$
    ie $\Delta''$ is a \tightening\ disc. It is also clear
    that $\bigl|\Delta''\cap\bigcup_i E_i\bigr|<\bigl|\Delta\cap\bigcup_i
    E_i\bigr|$, reducing complexity. We relabel $\Delta''$ as $\Delta$.
    \end{itemize}

In any case complexity is reduced, so eventually
$\Delta\cap\bigcup_i E_i=\partial\Delta$. Since we can regard the
$1$--handles of $\cal{H}_n$ as neighborhoods of the $E_i$'s, we
may assume that those do not intersect $\Delta$. In fact, because
of the exceptional 1--handle $e_l$ and disc $E_l$ which contain
$\partial\Delta$, this neighborhood argument yields $\Delta$
intersecting the union of the 1--handles of $\cal{H}_n$ only in a
collar neighborhood $F\subseteq\Delta$ of $\partial\Delta$. To
obtain the goal that $\Delta\cap\Lambda=\partial\Delta$ we will,
roughly, isotope $F$ along the $I$--fibers of the product
structure $D^2\times I$ of $e_l$ up to the point that
$\partial\Delta$ is contained in an extreme leaf of $\Lambda\cap
e_l$ (see \fullref{Rs:laminations}).

More precisely, consider an isotopy of $e_l$ which preserves both
of its product foliations $D^2\times I$ (preserving the
$I$--fibers) and takes $E_l$ to an extreme boundary leaf $E'_l$ of
$\Lambda\cap e_l$ (an isotopy along the $I$--fibers). Here we
assume that such a $E'_l$ is not any of the discs
$D^2\times(\partial I)\subseteq e_l$, possibly extending the
product structure of $e_l$ to a slightly larger neighborhood in
$H_n$. There are two choices of such a boundary leaf $E'_l$. We
choose the one that is ``to the side of $\Delta$'', in the sense
that $F$ is contained in the product $D^2\times [a,b]\subseteq
D^2\times I=e_l$ between $E_l$ and $E'_l$. The desired isotopy can
be obtained from an isotopy between $\text{Id}_I$ to a
homeomorphism $I\to I$ taking $a\mapsto b$ (or $b\mapsto a$,
depending on whether $E_l$ comes before of after $E'_l$ in the
orientation of $I$).

This isotopy fixes $D^2\times(\partial I)\subseteq e_l$, so it
extends to an isotopy of $H_n$ with support in $e_l$. We apply it
to the disc $\Delta$, obtaining $\Delta'$. By construction $F$ is
moved away from $\Lambda\cap e_l$, yielding
$\Delta'\cap\Lambda=\partial\Delta'$. Because $\Omega\cap e_l$
consists of $I$--fibers of $e_l$, which are preserved by the
isotopy, it also holds that $\nu(\Delta')=\nu(\Delta)$ and
$\Delta'$ is a tightening disc satisfying property 3. We relabel
$\Delta'$ as $\Delta$.

\textbf{Part 2: property \fullref{D:strong-p1}}\qua
Again let $n$ be such that $\Delta\subseteq H_n$. Let $L$ be the
leaf of $\Lambda\cap H_n$ containing $\partial\Delta$. Hence
$\partial\Delta$ bounds a disc $\Delta'\subseteq L$. Consider the
annulus $A=\overline{L-\Delta'}$ (\fullref{F:leaf}).

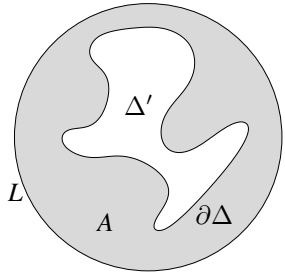
\begin{figure}[ht!]
\labellist\small
%\psfrag{D}{$\Delta'$}
%\psfrag{A}{$A$}
%\psfrag{L}{$L$}
%\psfrag{pD}{$\partial\Delta$}
\hair1pt
\pinlabel {$\Delta'$} at 143 180
\pinlabel {$A$} at 107 53
\pinlabel {$L$} [r] at 21 85
\pinlabel {$\partial\Delta$} [tl] at 202 73
\endlabellist
\centerline{\includegraphics[scale=0.35]{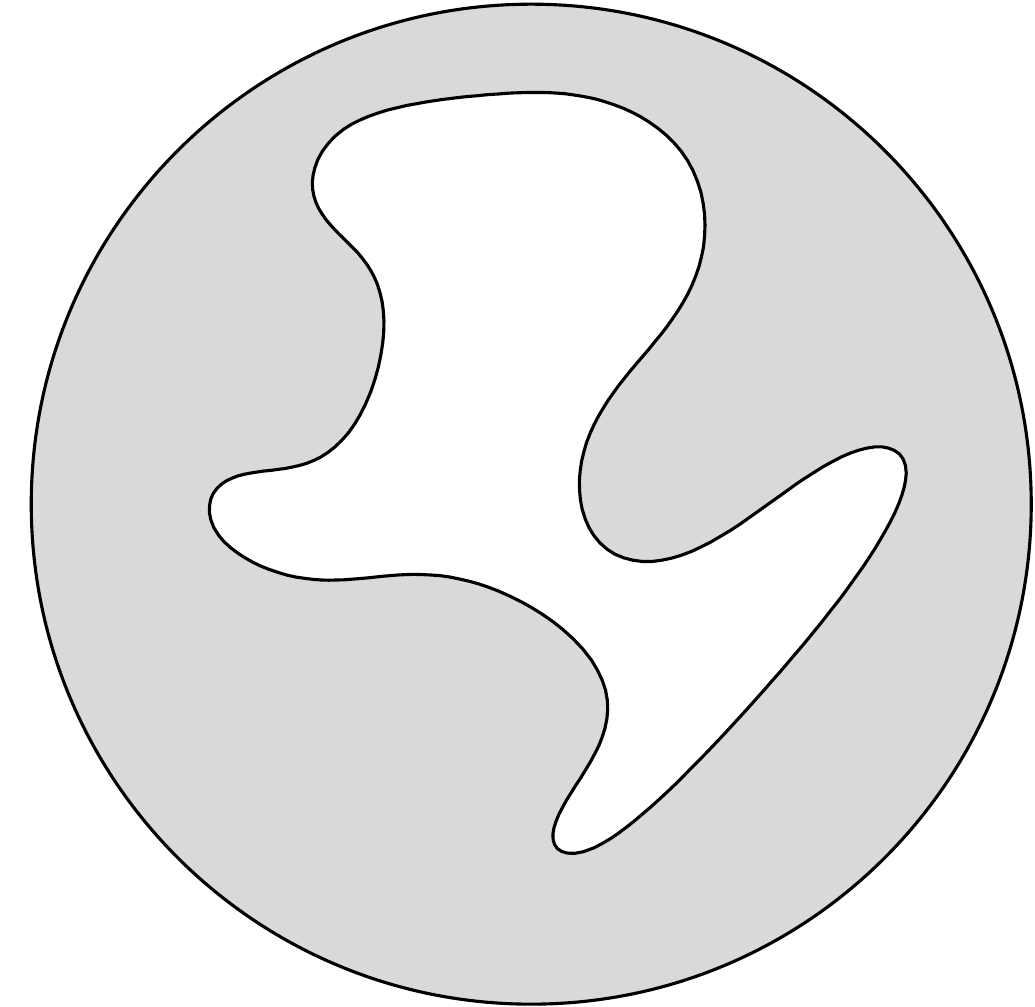}}
\caption{The curve $\partial\Delta$ separates a disc $\Delta'$ from an
annulus $A$ in $L$.}
\label{F:leaf}
\end{figure}

The desired disc $D$ (satisfying property 1) will be, essentially,
$\Delta\cup A$. Just push $\Delta\cup A$ slightly away from $L$,
leaving $\partial(\Delta\cup A)$ unchanged. Since $\Delta\cup A$
is transverse to $\Omega$ we can do that preserving
$\nu(\Delta\cup A)$. Therefore the resulting disc $D$ has the
property that $\nu(D)=\nu(\Delta\cup A)$.

It is easy to see that $D$ still is a tightening disc because the
whole operation increased intersection with $(\Omega,\nu)$ by the
same amount $\nu(A)$ both on the side of the disc and that of the
leaf. More precisely,
$$\nu(D)=\nu(\Delta)+\nu(A)<\nu(\Delta')+\nu(A)=\nu(L).$$
We relabel $D$ as $\Delta$, realizing property \fullref{D:strong-p1}
(and preserving property \fullref{D:strong-p3}).

\textbf{Part 3: property \fullref{D:strong-p4}}\qua
At this point we have a \tightening\ disc $\Delta$ and $n$ such
that $(\Delta,\partial\Delta)\subseteq(H_n,\partial H_n)$ and
$\Delta\cap\Lambda=\partial\Delta$.

It is clear that for a sufficiently large $m$, $\Delta$ does not
intersect the $0$--handles of $H_{-m}$, otherwise $\Delta$ would
intersect the singular set $S(\Omega)$ (see \fullref{Rs:laminations}). In fact, for a sufficiently large $m$,
$\Delta$ is transverse to the $I$--fibers of the $1$--handles
$D^2\times I$ of $H_{-m}$, otherwise $\Delta$ would not be
transverse to $\Omega$ (recall that for a $1$--handle $e_i$,
$\Omega\cap e_i$ consists of $I$--fibers). By taking a
sufficiently large $N$, $f^N(\Delta)$ will intersect $H_0$ only in
its $1$--handles and transverse to its $I$--fibers. We relabel
$f^N(\Delta)$ as $\Delta$ and $n+N$ as $n$.

Now note that, since $\Delta$ is transverse to the fibers of the
$1$--handles $e_i$ of $H_0$, $\Delta\cap e_i$ consists of discs
parallel to dual discs. We can then isotope $\Delta$ so that
$\Delta\cap e_i$ consists of dual discs. This isotopy does not
change $\inter{\Delta}{\Omega}{\nu}$ and $\Delta$ still is a
\tightening\ disc. From \eqref{Eq:inter}
$\inter{\Delta}{\cal{H}_0}{\widehat\nu}=\inter{\Delta}{\Omega}{\nu}$
(see \fullref{R:standard}), and it follows that $\Delta$ satisfies
property 4 of a strong tightening disc.

Properties \fullref{D:strong-p3} and \fullref{D:strong-p1} of $\Delta$
were preserved, so it only remains to verify property
\fullref{D:strong-p2}.

\textbf{Part 4: property \fullref{D:strong-p2}}\qua
We assume the following technical lemma, which will be proved later.

\begin{lemma}\label{L:nice discs}
Let $(E,\partial E)\subseteq(H_m,\partial H_m)$ be an embedded
disc intersecting $H_0$ in discs dual to the $1$--handles. If
$E\cap\Lambda$ (which may be empty) is contained in $\partial E$
then there exists a disc $E''\subseteq H_m$ such that
\begin{itemize}
   \item $E''\cap\Lambda\subseteq\partial E''=\partial E$,

   \item $E''\cap H_{m-1}$ consists of discs essential in $H_{m-1}$,

   \item $E''\cap H_0$ consists of dual discs, and

   \item $\inter{E''}{\cal{H}_0}{\widehat\nu}\leq
          \inter{E}{\cal{H}_0}{\widehat\nu}$.
\end{itemize}
\end{lemma}

So suppose that for some $0\leq j\leq n-1$, $\Delta\cap H_j$ does
not consist of essential discs. Let $m-1$ be the greatest such
value. Therefore the components of $\Delta\cap H_m$ are essential
discs. Apply the lemma for each such component $E$, replacing
$E\subseteq\Delta$ by $E''$. Now $\Delta\cap H_j$ consists of
essential discs for all $m-1\leq j\leq n$ and we proceed by
induction. In the end $\Delta$ will satisfy property 2 of the
definition of strong tightening disc.

We note that the other assertions of the lemma ensure that the
other properties \fullref{D:strong-p1}, \fullref{D:strong-p3} and
\fullref{D:strong-p4} of $\Delta$ are preserved.

This completes the proof of the proposition up to the proof of
\fullref{L:nice discs}.
\end{proof}

\begin{corollary}
The lamination is \tight\ if and only if there exists no strong
tightening disc.
\end{corollary}

\begin{proof}[Proof of \fullref{L:nice discs}]

Perturb $E$ in a neighborhood of $\partial H_{m-1}$ so that it is
transverse to $\partial H_{m-1}$. We want to change $E$ so that
$E\cap H_{m-1}$, which is a planar surface, consists just of
essential discs. A natural strategy would be the following:
isotope any inessential component of $E\cap H_{m-1}$ away from
$H_{m-1}$, so we would have to deal only with essential
components. These would be either discs, which we want, or
compressible in $H_{m-1}$. Simplify the compressible components by
compressing them. The problem with this last step is that these
compressions, being performed in $H_{m-1}$, could introduce
undesired intersections with $H_0$, increasing
$\inter{E}{\cal{H}_0}{\widehat\nu}$. Instead, we will perform surgeries
in $H_m-\ring H_{m-1}$, where they will not be subject to this
problem.

We start by making some quite general comments. ``Cut $H_m$ open
along $\Lambda$'' (ie take the completion of $H_m-\Lambda$
under a path-metric inherited from a metric in $H_m$) and let $C$
be the component that contains $E$. We claim that $\partial
H_{m-1}\cap C$ is incompressible in $C-\ring H_{m-1}$. Indeed, let
$F$ be a component of $\partial H_{m-1}\cap C$ and consider the
following commutative diagram:
$$\begin{CD}
\pi_1(F)   @>i_1>>    \pi_1(C-\ring H_{m-1})\\
@V i_2 V 1{-}1 V                                    @VV i_3V\\
\pi_1(\partial H_{m-1})         @>i_4>1{-}1>    \pi_1(H_m-\ring
H_{m-1})
    \end{CD}$$
where all homomorphisms are induced by the corresponding
inclusion. We claim that $i_2$ is injective. Indeed, recall that
$\Lambda\cap H_{m-1}$ consists os discs dual to the 1--handles in
$H_{m-1}$ (see \fullref{T:generic} and \fullref{Rs:laminations}) and hence $\Lambda\cap \partial H_{m-1}$
consists of essential curves in $\partial H_{m-1}$. Also $i_4$ is
injective (because $H_m-\ring H_{m-1}\simeq
\partial H_{m-1}\times I$). Therefore $i_4\circ i_2$, and thus
$i_3\circ i_1$, is injective. Then $i_1$ is injective and hence
$F$ is incompressible in $C-\ring H_{m-1}$.

From $\partial E\subseteq\partial H_m$ follows that
$E\cap(\partial H_{m-1}\cap C)$ consists of simple closed curves.
Therefore, if $S=E\cap(C-\ring H_{m-1})$, then $\partial S$ has
exactly one component in $\partial H_m\cap C$ (namely, $\partial
E$) and all the others in $\partial H_{m-1}\cap C$.

We now describe the operations that change $E$ to yield the
desired final $E''$. They will be performed in the interior of
$H_m\cap C$, therefore will not introduce intersections with
$\Lambda$ or with $\partial H_m$. For simplicity of notation, we
also label the modified discs by $E$.

\begin{process}
  Let $\gamma$ be a component of $E\cap(\partial H_{m-1}\cap
  C)$. If $\gamma\subseteq(\partial H_{m-1}\cap C)$ is not essential
  then it bounds a disc $ D\subseteq(\partial H_{m-1}\cap C)$. But
  $\gamma$ bounds a disc $ D'\subseteq E$ and $H_m$ is irreducible,
  therefore we can isotope $D'$ to $D$ (also pushing any part of $E$
  that is on the way).
\end{process}

Since $ D\cap H_0=\emptyset$ this operation does not introduce
intersections with $H_0$, hence it does not increase
$\inter{E}{\cal{H}_0}{\widehat\nu}$. It is clear that it also reduces
$|E\cap (\partial H_{m-1}\cap C)|$. By repeating the process we
can assume that $E\cap(\partial H_{m-1}\cap C)$ consists only of
curves that are essential in $\partial H_{m-1}\cap C$.

Recall $S=E\cap(C-\ring H_{m-1})$. Consider a component $F$ of
$S$. If $F$ is a disc then, by the previous paragraph, $\partial
F\subseteq (\partial H_{m-1}\cap C)$ is an essential curve. Hence
$F$ is a compressing disc for $\partial H_{m-1}\cap C$ in $C-\ring
H_{m-1}$, a contradiction. So $S$ cannot contain discs.

Now let $\gamma$ be a component of $E\cap(\partial H_{m-1}\cap C)$
which is innermost in $E$. It bounds a disc $E'\subseteq E$. From
the paragraph above $E'\subseteq(H_{m-1}\cap C)$, which will be
essential (by previous considerations). We conclude that when all
components of $E\cap(\partial H_{m-1}\cap C)$ are innermost in $E$
(in other words, when $S$ is connected) the proof is complete.

Suppose that that is not the case. The following procedure allows
us to assume that $S$ is incompressible (in $C-\ring H_{m-1}$).

\begin{process}
    If there exists a compressing disc $D$ for $S$ in
    $C-\ring H_{m-1}$ we ``compress'' along it in the following
    sense. Consider the disc $D'\subseteq E$ with
    $\partial D'=\partial D$. Irreducibility of $H_m$ implies that the
    sphere $D\cup D'$ bounds a ball. Change $E$ through an isotopy
    taking $D'$ to $D$ along this ball.
\end{process}

The isotopy does not increase $\inter{E}{\cal{H}_0}{\widehat\nu}$
because the compressing disc $D$ does not intersect $H_0$. Also,
note that $S\subseteq E$ is a planar surface therefore $\partial
D\subseteq S$ is separating, one of the sides being contained in
$D'\subseteq E$. Therefore the process, which essentially replaces
$D'$ with $D$, reduces $|\partial S|=|E \cap(\partial H_{m-1}\cap
C)|$. We can repeat the process and assume that $S$ is
incompressible in $C-\ring H_{m-1}$.

The next step makes $S$ connected and hence, as was argued before,
$E\cap H_{m-1}$ consisting of essential discs.

\begin{process}
  Supposing that $S$ is not connected, let $S'$ be a component such
  that $\partial S'\subseteq(\partial H_{m-1}\cap C)$. Now consider
  $S'$ in $H_m-\ring H_{m-1}$. We verify that $S'$, being
  incompressible in $C-\ring H_{m-1}$, is incompressible in
  $H_m-\ring H_{m-1}$. Indeed, if $ D$ is a compressing disc for $S'$ (in
  $H_m-\ring H_{m-1}$) we can simplify $ D\cap\Lambda$ by standard
  ``cut and paste'' techniques until $ D\cap\Lambda=\emptyset$.
  These surgeries will not change $\partial D$,
  so $ D$ still is a compressing disc. But now $ D\subseteq C$,
  contradicting incompressibility of $S'$ in $C-\ring H_{m-1}$. So
  $S'$ is incompressible in $H_m-\ring H_{m-1}$.

  But $H_m-\ring H_{m-1}$ has a product structure $(H_m-\ring
  H_{m-1})\simeq\partial H_{m-1}\times I$. It is a well known fact
  that an incompressible surface $S'$ in such a product with $\partial
  S'\subseteq\partial H_{m-1}\times\{0\}$ is parallel to a surface
  $S''\subseteq\partial H_{m-1}$. Let $P$ be the product bounded
  by $S'\cup S''$. It is easy to see that
  $P\cap\Lambda=\emptyset$, so we can use the $I$--fibers of $P$ to
  isotope $E\cap P$
  vertically, through an isotopy taking $S'$ a bit further than
  $S''$.
\end{process}

We note that this operation reduces $|\, S\,|$, does not change
intersections with $H_0$ and does not introduce intersections of
$E$ with $\Lambda$ (recall that $P\cap\Lambda=\emptyset$).
Repetition of the process yields $|\, S\,|=1$, ie $S$ is
connected and, therefore, $E\cap H_{m-1}$ consists of essential
discs.

We relabel $E$ as $E''$. It satisfies all the desired conditions
in the statement of the lemma, completing the proof.
\end{proof}

%%%%%%%%%%%%%%%%%%%%%%%%%%%%%%%%%%%%%%%%%%%%%%%%%%%%%%%%%%%%%%%%%%%%%%%%%%%%%%%%%%%%%%%%%%
%%%%%%%%%%%%%%%%%%%%%%%%%%%%%%%%%%%%%%%%%%%%%%%%%%%%%%%%%%%%%%%%%%%%%%%%%%%%%%%%%%%%%%%%%%

\subsection{Minimality and tightness}\label{SS:tightness}

We will use strong tightening discs to perform isotopies that
reduce the growth rate.

\begin{lemma}\label{L:reduce lambda}
If there is a strong tightening disc $\Delta\subseteq H_1$ then
$\lambda$ is not minimal.
\end{lemma}

\begin{obs}
\fullref{E:loose} shows a particular case of the
proposition. The proof follows essentially the same argument that
we gave in the example.
\end{obs}

\begin{proof}
Let $(\cal{H}_0,\widehat\nu)$ be the standard weighted handle
decomposition of $H_0$ (see \fullref{D:standard}) and
consider the weighted decomposition $(\cal{H}_1,\lambda\widehat\nu)$.
Since $\partial \Delta\subseteq H_1\cap\Lambda$ then $\partial
\Delta=\partial E_{i_0}^1$, where $E_{i_0}^1=f(E_{i_0})$ is a disc
dual to a 1--handle of $H_1$. Moreover, from \eqref{Eq:inter} (see
\fullref{R:standard}):
$$\inter{\Delta}{\cal{H}_0}{\widehat\nu}<
  \inter{E_{i_0}^1}{\cal{H}_0}{\widehat\nu}=\lambda\widehat\nu_{i_0}.$$
Since $\Delta\cap H_0$ consists of discs dual to the 1--handles of
$\cal{H}_0$ we can change $f$ through an isotopy taking $E_{i_0}$
to $\Delta$, without changing the handle decomposition and
preserving ``compatibility'' (see \fullref{Rs:laminations}).

Let $M'=(m'_{ij})$ be the new incidence matrix. Then
$m'_{ij}=m_{ij}$ if $i\neq i_0$ and $m'_{i_0j}=m_{i_0j}+d_j$,
where $d_j\in\mathbb{Z}$ satisfies the following inequality:
$$\sum_{1\leq j\leq k}d_j\widehat\nu_j<0.$$
If $M'$ is irreducible then for $i\neq i_0$,
$$(M'\widehat\nu)_i=(M\widehat\nu)_i=\lambda\widehat\nu_i$$
but, for the row $i_0$,
$$(M'\widehat\nu)_{i_0}=(M\widehat\nu)_{i_0}+\sum_{1\leq j\leq k}
  d_j\widehat\nu_j<\lambda\widehat\nu_{i_0},$$
hence, by \fullref{P:subinvariance},
$$\lambda'=\lambda(M')<\lambda,$$
completing the proof in this case.

If $M'$ is reducible then, possibly permuting the indices, $M'$
may be written in the form
$$M'= \left(
\begin{tabular}{c|ccc}
$*$&&$*$&\\ \hline
&&&\\
$0$&&$M''$&\\
&&&\\
\end{tabular}
\right)$$
where, for some $1\leq r\leq k-1$, $M''$ is irreducible of
dimension $(k-r)$. From the original system of discs
$\cal{E}=\{E_1,\dots,E_k\}$ we pass to an irreducible subsystem
$\cal{E}''=\{E_{r+1},\dots,E_k\}$.  The transpose of the incidence
matrix for these discs is precisely $M''$. We will show that also
in this case $\lambda(M'')<\lambda(M)$.

Since $M$ is irreducible there exist $r+1\leq i_1\leq k$ and
$1\leq j_1\leq r$ such that $m_{i_1j_1}>0$. Let $\widehat{\widehat\nu}$ be
the $(k-r)$--vector defined by $\widehat{\widehat\nu}_i=\widehat\nu_{i+r}$,
$1\leq i\leq k-r$ (ie $\widehat{\widehat\nu}$ consists of the last
$k-r$ coordinates of $\widehat\nu$). It follows that, for $i\neq
i_0-r$, $i_1-r$,
$$\bigl(M''\,\widehat{\widehat\nu}\,\bigr)_i\leq\bigl(M\,\widehat\nu\,\bigr)_{i+r}=
\lambda\widehat\nu_{i+r}=\lambda\widehat{\widehat\nu}_i
$$
and, for $i=i_1-r$,
$$
\bigl(M''\,\widehat{\widehat\nu}\,\bigr)_{i_1-r}\leq(M\,\widehat\nu)_{i_1}
-\bigl(m_{(i_1j_1)}\bigr)\widehat\nu_{j_1}<
\bigl(M\,\widehat\nu\,\bigr)_{i_1}=\lambda\widehat\nu_{i_1}=\lambda\widehat{\widehat\nu}_{i_1-r}.
$$
If $1\leq i_0\leq r$ then it follows that
$\lambda''=\lambda(M'')<\lambda$, proving the lemma in this case.

If $r+1\leq i_0\leq k$, in addition to the inequalities above, we
further have (for $i=i_0-r$)
$$
\bigl(M''\,\widehat{\widehat\nu}\,\bigr)_{i_0-r}=\sum_{r+1\leq j\leq
k}(n_{i_0j}')\widehat\nu_j
         \leq \sum_{1\leq j\leq k}(n'_{i_0j})\widehat\nu_j
=\bigl(M'\,\widehat\nu\,\bigr)_{i_0}<\lambda\widehat\nu_{i_0}=\lambda\widehat{\widehat\nu}_{i_0-r},
$$
therefore $\lambda''<\lambda$, completing the proof of the lemma.
\end{proof}

The hypothesis of $\Delta$ being contained in $H_1$ in the
statement of \fullref{L:reduce lambda} is needed because the
growth-reducing move has to be performed equivariantly. The
difficulty in proving the conjecture is precisely in finding such
a disc in $H_1$. We can do that under the hypotheses of
\fullref{P:valence} (ie maximum $0$--handle valence
$3$), as stated in \fullref{L:valence} below. The reason is,
essentially, that under such conditions every properly embedded
essential disc in $H_0$ which is disjoint from the original disc
system is parallel to a disc of this system.

\begin{lemma}\label{L:valence}
Suppose that all $0$--handles of the handle decomposition
$\cal{H}_0$ have valence $2$ or $3$. If $\Lambda(\cal{H})$ is not
\tight\ then there exists a strong tightening disc in $H_1$.
\end{lemma}

\begin{proof}
Supposing that $\Lambda$ is not \tight\ we use \fullref{P:equivalence} to get a strong tightening disc. Among all
such discs let $\Delta$ have minimal \emph{height}, in the sense
that if $H_n$ contains a strong tightening disc then
$\Delta\subseteq H_n$. Let $n$ be the smallest integer with the
property that $\Delta\subseteq H_n$. The lemma states that $n=1$,
so assume otherwise that $n\geq 2$.

If such a disc $D$ does not intersect any $0$--handles then it is
contained in a $1$--handle of $H_{n-1}$. In this case it is clear
that $D$ is parallel to a dual disc. If $D$ intersects a
$0$--handle we can assume it is actually contained in it: using
the product structure on the $1$--handles that $D$ intersects, the
fact that $D\cap\Lambda=\emptyset$ and that $\Lambda\cap H_{n-1}$
contains a representative of any dual disc, we can isotope $D$
into that $0$--handle. Since the handle has valence at most $3$,
again $D$ is parallel to a disc dual to a $1$--handle. Therefore
$\Delta\cap H_{n-1}$ consists of discs parallel to the co-cores of
the $1$--handles.

Let $D$ be such a disc, parallel to the co-core
$E_i^{n-1}\subseteq H_{n-1}$. Since $D\subseteq\Delta$ (which is a
strong tightening disc) then $D\cap H_0$ consists of dual discs of
$\cal{H}_0$, so it makes sense to consider $\inter
{D}{\cal{H}_0}{\widehat\nu}$. If
\begin{equation}\label{Eq:valence}
\inter{D}{\cal{H}_0}{\widehat\nu}\geq
\inter{E_i^{n-1}}{\cal{H}_0}{\widehat\nu}
\end{equation}
for all discs $D\subseteq\Delta\cap H_{n-1}$ then we can alter
$\Delta$ by an isotopy in such a way that each $D$ is moved to the
corresponding dual disc. This operation does not increase
$\inter{\Delta}{\cal{H}_0}{\widehat\nu}$ and preserves the other
properties of strong tightening discs. Now
$\inter{E_i^{n-1}}{\cal{H}_0}{\widehat\nu}=\lambda^{n-1}\widehat\nu_i$ is
precisely the weight $\inter{E_i^{n-1}}{\Omega}{\nu}$ on the
$i$-th $1$--handle of $H_{n-1}$. If we apply $f^{-n+1}$ to
$\Delta$ and $E_i^{-n+1}$ the inequality
$$
\inter{f^{-n+1}(D)}{\cal{H}_0}{\widehat\nu}\geq\inter{E^i}{\cal{H}_0}{\widehat\nu},
$$
is obtained from \eqref{Eq:valence} by multiplying both sides by a
factor of $\lambda^{-n+1}$. That proves that
$f^{-n+1}(\Delta)\subseteq H_1$ is a strong tightening disc, a
contradiction to the assumption $n\geq 2$.

The argument above then shows that
$$
\inter{D}{\cal{H}_0}{\widehat\nu}<\inter{E_i^{n-1}}{\cal{H}_0}{\widehat\nu}
$$
for some $D\subseteq \Delta\cap H_{n-1}$.  Modifying $D$ through
an isotopy supported in a regular neighborhood of $\partial
H_{n-1}$ (hence preserving $\inter{D}{\cal{H}_0}{\widehat\nu}$), we
can assume that $\partial D=\partial E_i^{n-1}$. Now $D\subseteq
H_{n-1}$ is a strong tightening disc, contradicting minimality of
$n$.

Therefore $n=1$ and the proof is complete.
\end{proof}

We recall and prove:

\begin{prop-valence}
Suppose that all $0$--handles of the handle decomposition
$\cal{H}$ have valence $2$ or $3$. If $\Lambda(\cal{H})$ is not
\tight\ then $\lambda(\cal{H})$ is not minimal.
\end{prop-valence}
\begin{proof}
Use \fullref{L:valence} and apply \fullref{L:reduce lambda}.
\end{proof}

\begin{corollary}\label{C:genus 2}
Let $f\co H\to H$ be a generic automorphism of a handlebody of
genus $2$. If $\lambda$ is minimal then $\Lambda$ is \tight.
\end{corollary}
\begin{proof}
Consider the handle decomposition $\cal{H}_0$ of $H_0$ with
co-core $\cal{E}_0$. If every 0--handle has valence $2$ or $3$
then \fullref{P:valence} completes the proof. So we assume
that this is not the case. Since $H$ has genus $2$, then it has
just one 0--handle with valence $4$ and the others (possibly none)
with valence $2$.

We sketch the proof in the one--dimensional setting: we consider
$\Gamma_0$ the graph corresponding to $\cal{H}_0$, which will have
a valence $4$ vertex and some valence $2$ vertices. Using the
height function in $H_1-\ring{H}_0$ (projection on the $I$
coordinate of the product) we obtain graphs $\Gamma_t$ dual to
$\cal{E}_1\cap H_t$ (when $t$ is a regular value of the height
function in $\cal{E}_1$). As $t$ increases, $\Gamma_t$ changes by
folds. The first fold will have to happen at the valence $4$
vertex (folds at valence $2$ vertices do not happen by
incompressibility). We want to say that this fold replaces the
valence $4$ vertex by two of valence $3$, reducing the problem to
the previous case. It could happen that the fold is done along an
edge determining a closed loop and, after that, the fold actually
replaces the vertex by another one still with valence $4$. We can
solve this problem by sufficiently subdividing the edges of
$\Gamma$ (ie introducing valence $2$ vertices in the interior
of the edges).

We now give the more detailed proof in the two--dimensional
setting. The fact that one $0$--handle has valence $4$ means, in
genus $2$, that $\mathcal{E}_0$ consists of just two isotopy
classes of discs. We want each of these isotopy classes to contain
at least two distinct discs of $\mathcal{E}_0$. This can be easily
achieved by \emph{splittings} (see the paper \cite{UO:Autos} by Oertel,
or the dual \emph{subdivision} of Bestvina and Handel \cite{BH:Tracks}).

Since $\lambda$ is minimal, $\Lambda$ has the incompressibility
property. In particular $\cal{E}_1-\ring{H}_0$ is incompressible
in $H_1-\ring{H}_0$. We then have a height function in
$H_1-\ring{H}_0$ with respect to which we may suppose that
$\cal{E}_1$ is in Morse position having just saddles as critical
points, no pair in the same level (see \cite{UO:Autos}). Let $t$
be a bit greater than the first critical value and consider the
intermediate $H_0\subsetneq H_t\subsetneq H_1$. Then
$\cal{E}_t=\cal{E}_1\cap H_t$ is a system of discs for $H_t$. It
determines a handle decomposition $\cal{H}_t$. By regarding $H_t$
as a ``new $H_0$'' and following the construction of the invariant
measures laminations, we obtain the same original lamination and
growth: $(\Lambda_t,\mu_t)=(\Lambda,\mu)$,
$(\Omega_t,\nu_t)=(\Omega,\nu)$ and $\lambda_t=\lambda$ (see
\cite{UO:Autos} for details).

We claim that $\cal{E}_t$ contains three distinct isotopy classes
of discs. To see this we identify $H_t$ with $H_0$ through the
product structure and note that $\cal{E}_t$ is obtained from
$\cal{E}_0\subseteq H_0$ by replacing two distinct discs by their
band sum. This sum is done along a band contained in a
$0$--handle. If such a $0$--handle had valence $2$ then the band
sum would join two parallel discs, contradicting incompressibility
in $H_1-\ring{H}_0$, so the band is contained in the $0$--handle
with valence $4$. But the band sum of discs in a $0$--handle with
valence greater then $3$ yields a disc in a new isotopy class. Now
recall that we chose $\cal{E}_0$ to have at least two discs in
each class. Therefore a single band sum will preserve at least one
disc in each of these original classes, proving the claim that
$\mathcal{E}_t$ contains three classes of discs.

Now that $\cal{E}_t$ contains three distinct isotopy classes of
discs, all its $0$--handles have valence $3$ or $2$. Therefore, by
\fullref{P:valence}, $\Lambda_t=\Lambda$ is tight.
\end{proof}
\begin{obs}
There is another interesting point of view, from which we also
sketch the proof here. It uses the disc complex of $H$. The disc
complex $\cal D(H)$ of a handlebody $H$ of genus $2$ has dimension
$2$. We identify each $H_t$ with $H$ by collapsing the $I$--fibers
of the product $H-\mathring{H}_t\simeq \partial H\times I$. Now
consider the path $t\mapsto \Lambda\cap H_t$ on $\cal D(H)$ (here
we use normalized transverse measures as barycentric coordinates).
But $\Lambda\cap H_t$ determines a complete system of discs for
any $t$ so it has at least two isotopy classes of discs. Therefore
the path $t\mapsto \Lambda\cap H_t$ never intersects the
$0$--skeleton of $\cal D(H)$. But the path intersects infinitely
many simplices, so it has to intersect the interior of some
$2$--simplex of $\cal D(H)$, which corresponds to three distinct
isotopy classes of discs. The argument is finished as before.
\end{obs}

%%%%%%%%%%%%%%%%%%%%%%%%%%%%%%%%%%%%%%%%%%%%%%%%%%%%%%%%%%%%%%%%%%%%%%%%%%%%%%%%%%%%%%%%%%
%%%%%%%%%%%%%%%%%%%%%%%%%%%%%%%%%%%%%%%%%%%%%%%%%%%%%%%%%%%%%%%%%%%%%%%%%%%%%%%%%%%%%%%%%%

\subsection{Tightness implies minimality}\label{SS:minimality}

The next goal is to prove \fullref{T:minimality}. For that, we
need a technical result coming below. If $G$, $G'\subseteq H$ are
embedded graphs we shall say that $G'$ \emph{follows $G$} if
$G'\subseteq F_G$, where $F_G$ is a fixed fibered neighborhood of
$G$, the vertices of $G'$ are contained in the union of the
neighborhoods of vertices of $G$ and the edges of $G'$ are
transverse to the fibers of $F_G$ over the edges.

Now suppose that $g\co H\to H$ is an automorphism and that $g(G)$
follows $G$. In this case we define an incidence matrix $N=N_G(g)$
by $n_{ij}=|g(e_i)\cap F(e_j)|$, where $F(e_j)\subseteq F_G$ is
the fibered neighborhood over the edge $e_j$. If $N$ is
irreducible we say that its Perron--Frobenius eigenvalue
$\lambda_G=\lambda_G(g)$ is \emph{the growth of $g$ on $G$}.

In the following we assume that $f\co H\to H$ and handle
decomposition $\cal{H}_0$ of $H_0$ are fixed. We consider the
corresponding disc system $\cal{E}_0$, dual graph $\Gamma_0$,
laminations $(\Lambda,\mu)$, $(\Omega,\nu)$ and growth rate
$\lambda$.

\begin{proposition}\label{P:following}
Suppose that there exists $g$ isotopic to $f^{-1}$ and graph
$G\subseteq H$ such that $g(G)$ follows $G$ with
$\lambda_G<\lambda$. If $\Gamma_0$ is isotopic to a graph
$\Gamma_0'$ which follows $G$ then $\Lambda$ is not tight.
\end{proposition}
\begin{proof}
As usual, the measured lamination $(\Omega,\nu)$ determines the
standard weighted handle decomposition $(\cal{H}_0,\widehat\nu)$ of
$H_0$. We use $(\cal{H}_0,\widehat\nu)$ and $f^{-n}$ to induce a
weighted handle decomposition $(\cal{H}_{-n},\widehat\nu)$ in
$H_{-n}=f^{-n}(H_0)$ (ie the weight in the $1$--handle
$f^{-n}(e_i)$ of $H_{-n}$ is $\widehat\nu_i$).

Let $E_i$ be the co-core of the $1$--handle $e_i$ of $\cal{H}_0$.
By the eigenvalue property of $\widehat\nu$, for any $n>0$
$$
\frac{\inter{E_i}{\cal{H}_{-n}}{\widehat\nu}}{\lambda^n}=\widehat\nu_i.
$$

The goal of the argument is to find a disc $\Delta$ isotopic to
$E_i$ (rel $\partial E_i$) such that, for some $N>0$,
\begin{equation}\label{Eq:goal}
\frac{\inter{\Delta}{\cal{H}_{-N}}{\widehat\nu}}{\lambda^N}<\widehat\nu_i.
\end{equation}
Since $E_i$ may be chosen as a leaf of $\Lambda\cap H_0$ and
$
\inter{\Delta}{\Omega}{\nu}=\frac{\inter{\Delta}{\cal{H}_{-N}}{\widehat\nu}}{\lambda^N}$,
such a $\Delta$ will then be a tightening disc.

The rough strategy is to isotope $H_0$ into the fibered
neighborhood $F_G$ of $G$ and iterate $g$. Since its growth is
smaller, the number of components of the intersection with the
disc $E_i$ will grow more slowly than originally, what will yield
a tightening disc. We shall develop this idea more precisely.

A big part of the proof consists of certain constructions, as
follows. We choose an isotopy $h$ taking $\Gamma_0$ to
$\Gamma_0'$. We now consider $H_0$ as a neighborhood of $\Gamma_0$
and can obtain an isotopic $H_0'\subseteq F_G$. The weighted
handle decomposition of $H_0$ determines a handle decomposition
$(\cal{H}_0',\widehat\nu)$ of $H_0'$. Also, the fibered structure of
$F_G$ determines a handle decomposition $\cal{G}$ of $F_G$ in the
natural way (ie neighborhoods of vertices correspond to
$0$--handles and fibered neighborhoods over the edges to
$1$--handles). By adjusting $h$ we can assume further that
$\cal{H}_0'$ and $\cal{G}$ are \emph{compatible} in the following
sense. Any dual disc of a $1$--handle of $\cal{G}$ intersects
$H_0'$ in dual discs of $\cal{H}_0'$.

Now use $(\cal{H}_0',\widehat\nu)$ to induce a system of weights on
$\cal{G}$ in the following way. Let $e_0,\dots, e_l$ be the
$1$--handles of $\cal{G}$ and, for each $i$, let $D_i$ be a dual
disc of $e_i$. Define
\begin{equation}\label{Eq:new weights}
(\widehat v_G)_i=\inter{D_i'}{\cal{H}_0'}{\widehat\nu}.
\end{equation}
Such a $\widehat v_G$ is well defined: the way the handles of
$\cal{G}$ intersect those of $\cal{H}_0'$ assures that
\eqref{Eq:new weights} above makes sense and does not depend on
the choice of disc dual to $e_i$. This defines a weighted
decomposition $(\cal{G},\widehat v_G)$.

We recall $g$ from the hypotheses of the lemma. The decomposition
$\cal{G}$ determines a decomposition $\cal{G}_{-1}$ of
$g\left(F_G\right)$. We can adjust $g$ so that $\cal{G}$ and
$\cal{G}_{-1}$ are compatible (in particular,
$g\left(F_G\right)\subseteq F_G$). Through $g^n$ we can define
decompositions $\cal{G}_{-n}$ of $g^n\left(F_G\right)$. It is
clear that these decompositions are automatically compatible.

We now define $H_{-n}'=g^n(H_0')$, also with weighted handle
structure $(\cal{H}_{-n}',\widehat\nu)$ defined through $g^n$. Clearly
$(\cal{H}_{-n}',\widehat\nu)$ is isotopic to $(\cal{H}_{-n},\widehat\nu)$.
It is also clear that if $-m\leq -n$ then $\cal{G}_{-n}$ and
$\cal{H}_{-m}'$ are compatible.

Recall the weight system $\widehat v_G$ defined in \eqref{Eq:new
weights} above. Now define weighted handle decompositions
$(\cal{G}_{-n},\widehat v_G)$. For a general surface $S$ the
construction implies that (assuming that both sides makes sense)
\begin{equation}\label{Eq:induced}
\inter{S}{\cal{H}_{-n}'}{\widehat\nu}=\inter{S}{\cal{G}_{-n}}{\widehat
v_G}.
\end{equation}
To complete the constructions necessary in this proof, let $E_i$
be a disc in the original system $\cal{E}_0$. For a technical
reason we extend it through the product structure in
$H-\ring{H}_0$ to a disc $(E,\partial E)\subseteq (H,\partial H)$.
It is clear that
$\inter{E}{\Omega}{\nu}=\inter{E_i}{\Omega}{\nu}=\widehat\nu_i$. We
isotope $E$ and also assume that $E\cap F_G$ consists of dual
discs of $\cal{G}$.

Recall that $\lambda_G=\lambda_G(g)$ is the growth of $g$ on $G$.
Clearly the sequence
$$n\mapsto \frac{\inter{E}{\cal{G}_{-n}}{v_G}}{(\lambda_G)^n}$$
is bounded. But $\lambda>\lambda_G$, therefore
$$\frac{\inter{E}{\cal{G}_{-n}}{v_G}}{\lambda^n}\to 0$$
and then, for some $N$,
\begin{equation}\label{Eq:got small}
\frac{\inter{E}{\cal{G}_{-N}}{v_G}}{\lambda^N}<\widehat\nu_i.
\end{equation}
By \eqref{Eq:induced} and \eqref{Eq:got small},
\begin{equation}\label{Eq:almost tightening}
\frac{\inter{E}{\cal{H}_{-N}'}{\widehat\nu}}{\lambda^N}<\widehat\nu_i.
\end{equation}
Now note that there is an ambient isotopy $h'\co H\to H$ such that
$h'(\cal{H}_{-N}',\widehat\nu)=(\cal{H}_{-N},\widehat\nu)$. By applying
$h'$ to $(\cal{H}_{-N}',\widehat\nu)$, it follows from \eqref{Eq:almost
tightening} above that
\begin{equation}\label{Eq:tightening}
\frac{\inter{h'(E)}{\cal{H}_{-N}}{\widehat\nu}}{\lambda^N}<\widehat\nu_i.
\end{equation}
We can choose $h'$ restricting to the identity at $\partial H$, so
$\Delta'=h'(E)$ has the property that $\partial\Delta'=\partial
E$.

Now recall that $E$ is the extension to $H$ of the co-core $E_i$
of a handle of $\cal{H}_0$. We can use the product structure in
$H-\ring{H}_0$ to obtain a disc $\Delta\subseteq H_0$ from
$\Delta'$. Clearly $\partial\Delta=\partial E_i$, $E_i$ may be
chosen as a leaf of $\Lambda\cap H_0$ and
$\inter{\Delta}{\cal{H}_{-N}}{\widehat\nu}=\inter{\Delta'}{\cal{H}_{-N}}{\widehat\nu}$.
Therefore \eqref{Eq:goal} follows from \eqref{Eq:tightening},
showing $\Delta$ as a \tightening\ disc.
\end{proof}

\begin{thm-minim}
If $\Lambda$ is tight then $\lambda$ is minimal.
\end{thm-minim}
\begin{proof}
It is a corollary of \fullref{P:following} above. We prove
the countrapositive, so assume that $\lambda$ is not minimal. Then
there exists another structure $\widehat{\cal{H}}_0$ for some $\widehat
H_0$ and representative $\widehat f$ for which the growth rate $\widehat
\lambda$ is less than $\lambda$. We consider the graph $\widehat
\Gamma_0$ corresponding to $\widehat{\cal{H}}_0$. It is direct that 1)
$(\widehat f)^{-1}(\widehat \Gamma_0)$ follows $\Gamma_0$ (with
$\lambda_{\widehat \Gamma_0}=\widehat \lambda<\lambda$) and 2) that
$\Gamma_0$ is isotopic to a $\Gamma_0'=\widehat \Gamma_0$. By
\fullref{P:following} $\Lambda$ is not tight, completing
the proof.
\end{proof}

\begin{obs}
\fullref{T:minimality} can be used to find the minimal growth
of some actual examples (see the author's doctoral thesis \cite{LC:thesis}).
\end{obs}

%%%%%%%%%%%%%%%%%%%%%%%%%%%%%%%%%%%%%%%%%%%%%%%%%%%%%%%%%%%%%%%%%%%%%%%%%%%%%%%%%%%%%%%%%%
%%%%%%%%%%%%%%%%%%%%%%%%%%%%%%%%%%%%%%%%%%%%%%%%%%%%%%%%%%%%%%%%%%%%%%%%%%%%%%%%%%%%%%%%%%

\subsection {Applications: comparing growth rates}\label{SS:boundary}

In this subsection we prove some results on the growth rate of an
irreducible automorphism with respect to a \tight\ lamination.

The following is a direct corollary of \fullref{T:minimality}.

\begin{corollary}\label{C:growth power}
If $\Lambda$ is \tight\ and $n>0$ then
$\lambda_{min}(f^n)=\lambda^n$.
\end{corollary}
\begin{proof}
Since $f^n(\Lambda,\mu)=(\Lambda,\lambda^n\mu)$ and
$f^n(\Omega,\nu)=(\Omega,\lambda^{-n}\nu)$, the growth of $f^n$
with respect to the handle decomposition that determine $\Lambda$
and $\Omega$ is $\lambda^n$. But $\Lambda$ is \tight, therefore
$\lambda^n$ is minimal.
\end{proof}

In Oertel's paper \cite{UO:Autos} the question is posed whether there is
any relation between the growth rate $\lambda=\lambda(f)$ of a generic
automorphism and the growth rate $\lambda_\partial$ of the pseudo-Anosov
restriction $\partial f=f|_{\partial H}$ to the boundary. The following
is another corollary of \fullref{P:following}.

\begin{corollary}\label{C:small growth}
Let $f\co H\to H$ be a generic automorphism. If $\Lambda$ is tight
then $\lambda\leq\lambda_\partial$.
\end{corollary}
\begin{proof}
We prove the countrapositive, so suppose that
$\lambda_\partial<\lambda$.

Let $g=f^{-1}$ and let $\tau$ be a stable train-track of $\partial
g$ (we regard $\tau$ simply as a graph, with switches for vertices
and branches for edges). Let $F_\tau\subseteq H$ be a fibered
neighborhood of $\tau$ and isotope $g$ so that $g(\tau)$ follows
$\tau$. It is clear that $\Gamma_0$ is isotopic to a graph
$\Gamma_0'$ which follows $\tau$. Indeed, $\Gamma_0$ is boundary
parallel and $\tau$ fills $\partial H$.
\fullref{P:following} completes the proof.
\end{proof}

\begin{corollary}\label{C:growth genus2}
The (minimal) growth rate of a generic automorphism $f\co H\to H$
of a handlebody of genus $2$ is less then or equal to the growth
rate of the pseudo-Anosov $\partial f=f\vert_{\partial H}$.
\end{corollary}
\begin{proof}
Use \fullref{C:small growth} and \fullref{C:genus 2}.
\end{proof}

%%%%%%%%%%%%%%%%%%%%%%%%%%%%%%%%%%%%%%%%%%%%%%%%%%%%%%%%%%%%%%%%%%%%%%%%%%%%%%%%%%%%%%%%%%
%%%%%%%%%%%%%%%%%%%%%%%%%%%%%%%%%%%%%%%%%%%%%%%%%%%%%%%%%%%%%%%%%%%%%%%%%%%%%%%%%%%%%%%%%%

\bibliographystyle{gtart}
\bibliography{link}

\end{document}